\documentclass[10pt,a4paper]{amsart}
\usepackage[UKenglish]{babel}
\usepackage{amssymb}
\usepackage{mathrsfs} 
\usepackage[T1]{fontenc}
\usepackage{ae,aecompl}
\usepackage{times}

\usepackage{tikz}

\usetikzlibrary[arrows]
\usetikzlibrary{decorations.pathmorphing}
\usetikzlibrary{decorations.shapes}
\usetikzlibrary{shapes,snakes}

\usepackage{hyperref}
\usepackage[alphabetic,backrefs]{amsrefs}
\usepackage{enumerate}
\usepackage[cmtip,all]{xy}
\usepackage{booktabs}
\usepackage[top=2cm,bottom=2cm,left=2.5cm,right=2.5cm]{geometry}

\newcommand{\rd}{{\rm d}}

\newcommand{\rD}{{\rm D}}

\newcommand{\rG}{{\rm G}}

\newcommand{\rU}{{\rm U}}

\newcommand{\bi}{{\bf i}}

\newcommand{\bA}{{\bf A}}

\newcommand{\bE}{{\bf E}}

\newcommand{\bP}{{\bf P}}

\newcommand{\cA}{\mathcal{A}}

\newcommand{\cE}{\mathcal{E}}

\newcommand{\cG}{\mathcal{G}}

\newcommand{\cM}{\mathcal{M}}

\newcommand{\cO}{\mathcal{O}}

\newcommand{\cS}{\mathcal{S}}

\newcommand{\cV}{\mathcal{V}}

\newcommand{\fg}{{\mathfrak g}}
\newcommand{\frg}{{\mathfrak g}}

\newcommand{\N}{\mathbb{N}}
\newcommand{\Z}{\mathbb{Z}}
\newcommand{\Q}{\mathbb{Q}}
\newcommand{\R}{\mathbb{R}}
\newcommand{\C}{\mathbb{C}}

\newcommand{\Spin}{{\rm Spin}}

\renewcommand{\P}{\bP}

\newcommand{\Aut}{\mathrm{Aut}}
\newcommand{\Crit}{\mathrm{Crit}}
\renewcommand{\det}{\mathop\mathrm{det}\nolimits}

\newcommand{\End}{{\mathrm{End}}}

\renewcommand{\epsilon}{\varepsilon}

\newcommand{\Hol}{\mathrm{Hol}}

\newcommand{\Lie}{\mathrm{Lie}}

\newcommand{\delbar}{\bar \del}
\newcommand{\del}{\partial}

\newcommand{\id}{\mathrm{id}}

\renewcommand{\Im}{\mathop{\mathrm{Im}}}

\newcommand{\loc}{{\rm loc}}

\renewcommand{\Re}{\mathop{\mathrm{Re}}}

\newcommand{\tr}{\mathop{\mathrm{tr}}\nolimits}

\newcommand{\qandq}{\quad\text{and}\quad}

\def\<{\mathopen{}\left<}
\def\>{\right>\mathclose{}}
\def\({\mathopen{}\left(}
\def\){\right)\mathclose{}}

\usepackage{multicol, color}

\definecolor{gold}{rgb}{0.85,.66,0}
\definecolor{cherry}{rgb}{0.9,.1,.2}
\definecolor{burgundy}{rgb}{0.8,.2,.2}
\definecolor{orangered}{rgb}{0.85,.3,0}
\definecolor{orange}{rgb}{0.85,.4,0}
\definecolor{olive}{rgb}{.45,.4,0}
\definecolor{lime}{rgb}{.6,.9,0}
\definecolor{green}{rgb}{.2,.7,0}
\definecolor{grey}{rgb}{.4,.4,.2}
\definecolor{brown}{rgb}{.4,.3,.1}

\newtheorem{theorem}{Theorem}
\newtheorem{proposition}{Proposition}
\newtheorem{corollary}[proposition]{Corollary}
\newtheorem{lemma}[proposition]{Lemma}

\theoremstyle{remark}
\newtheorem{remark}[equation]{Remark}

\theoremstyle{definition}
\newtheorem{definition}[proposition]{Definition}
\newtheorem{example}[proposition]{Example}

\def\tr{{\textrm{Tr}}}

\def\de{\rd}
\def\De{\rD}
\def\Dex{\De_\xi}

\newcommand{\feps}{f_\epsilon}
\newcommand{\xeps}{X_\epsilon}

\begin{document}

\title{Gauge theory and ${\rm G}_2$--geometry on Calabi-Yau  links}
\author{Omegar Calvo-Andrade, L\'azaro O. Rodr\'iguez D\'iaz, Henrique N. S\'a Earp}
\date{\today}
\subjclass[2010]{53C07, 53C25, 57R15; 32S55, 57R18}

\begin{abstract}
    The $7$--dimensional link $K$ of a weighted homogeneous hypersurface on the round $9$--sphere in $\C^5$ has a  nontrivial null Sasakian structure which is contact Calabi-Yau, in many cases. It admits a canonical co-calibrated $\rm G_2$--structure $\varphi$ induced by the  Calabi-Yau $3$--orbifold basic geometry. We distinguish these pairs $(K,\varphi)$ by the Crowley-Nordstr\"om $\Z_{48}$--valued $\nu$ invariant, for which we prove odd parity and provide an algorithmic formula.  

We describe  moreover a natural Yang-Mills theory on such spaces, with many  important features of the torsion-free case, such as a Chern-Simons formalism and topological energy bounds. In fact,  compatible  $\rm G_2$--instantons on holomorphic Sasakian bundles over $K$ are exactly the transversely  Hermitian Yang-Mills connections. As a proof of principle, we obtain  $\rm G_2$--instantons over the Fermat quintic link from stable bundles over the smooth projective Fermat
quintic, thus relating in a concrete example the Donaldson-Thomas theory of the quintic threefold with a conjectural $\rm G_2$--instanton  count.     
\end{abstract}

\maketitle

\tableofcontents
\newpage
\section{Introduction}
We propose a contemporary angle on Milnor's celebrated study of singular
hypersurface links  \cite{Milnor1969},
from the  perspective of special metrics and higher-dimensional
gauge theory. Our starting point is the observation that several topological properties of the Milnor fibre and its boundary the link (see Section \ref{sec-geometric-structures}) resemble those of the ${\rm G}_2$--invariant $\nu$ recently introduced by Crowley and Nordstr\"{o}m \cite{Crowley2015b}, suggesting to optimists that Milnor's construction might
 be related to   ${\rm G}_2$--geometry.

\subsection{$\rG_2$--metrics on Calabi-Yau links}

Let $\cV \subset\C^{n+1}$ be a complex analytic variety with an isolated singularity at the origin. Milnor proved that $\cV$ intersects transversally every sufficiently small sphere $S^{2n+1}:=\partial B_{\epsilon}(0)$, and
the \emph{link} 
$$
K:= \cV\cap S^{2n+1}
$$ is a $(n-2)$--connected smooth manifold with $\dim_\R K=2\dim_\C\cV-1$. The topologies of $\cV$ and of its embedding in $\C^{n+1}$ are completely determined by the embedding $K\hookrightarrow S^{2n+1}$.

Suppose henceforth that  $\cV=(f)$ is an affine hypersurface  defined by a 
homogeneous polynomial $f:\C^{n+1}\rightarrow\C$, with $f(0)=0$
and $\Crit(f)\cap B_{\epsilon}(0)=\{0\}$.  The Hopf fibration
$\pi:S^{2n+1}\rightarrow \P^n$
characterises the corresponding link  $K_f$ in a natural way as the total space of a $S^1$--bundle over the smooth projective hypersurface $V$ defined by $f$:
$$
\pi:K_f \stackrel{S^1}{\longrightarrow} V\subset\P^n.
$$ 
As a circle bundle, $K_f$ carries a global angular form $\theta\in\Omega^{1}(K)$, whose restriction to each fibre $\pi^{-1}(x)$ generates the cohomology $H^1(\pi^{-1}(x),\R)$. Its exterior derivative $\de\theta=-\pi^*e\in \Omega^2(K)$ is the pullback of the Euler class on the base (compare with Lemma \ref{lemma: dtheta is (1,1)}, below).

If the link has degree $n+1$, then the projective variety $V$ is a Calabi-Yau $(n-1)$--fold. Fixing $n=4$, a \emph{quintic} (possibly weighted) link $K_f $ is a smooth Sasakian $7$-manifold
fibering by circles over the smooth Calabi-Yau $3$--fold $V$, and it is the simplest example of a
 \emph{Calabi-Yau (CY) link} (see Definition \ref{def: CY link}).  Now, it is well-known  that the Riemannian product of a Calabi-Yau $3$--fold
and a circle carries a torsion-free $\rG_2$--structure, so we define naturally (see also Theorem 2.5 in \cite{Gray1969}) the following $\rG_2$-structure on $K_f $: 
\begin{equation}
\label{eq-GrayG2structure}
\begin{array}{rcl}        
        \varphi &:=& \theta\wedge \omega + \Im\epsilon,\\ 
        \psi &:=&\frac{1}{2}\omega\wedge\omega 
        + \theta\wedge\Re\epsilon=*\varphi
\end{array}
\end{equation}
where $\omega$ and $\epsilon$ are respectively the K\"{a}hler and holomorphic volume forms defining the Calabi-Yau structure on $V$ and we denote identically
differential forms and their pullbacks under $\pi$.
Although in the nontrivial fibration case this structure has torsion, it is still cocalibrated (see Section \ref{sec: G_2-geometry of links}), in the sense that $\de\psi=0$:
\begin{theorem}
\label{thm: cocalibrated G2-structure}
 Every quintic link $K_f $ is a  $2$--connected, compact, smooth real
 $7$--manifold admitting the natural cocalibrated $\rG_2$--structure 
(\ref{eq-GrayG2structure}).
\end{theorem}

It should be stressed that Theorem \ref{thm: cocalibrated G2-structure} has recently been found and subsumed, independently, by Habib and Vezzoni \cite[\S6.2]{Habib2015} in the context of \emph{contact Calabi-Yau} (cCY) geometry. Their theory extends the above discussion to weighted homogenous links and therefore yields many more examples of CY links, fibering over CY $3$--orbifolds in weighted $\P^4(w)$ (see Section \ref{sec: contact Calabi Yau}). This is very fortunate, because otherwise the Fermat quintic would the only strictly homogeneous quintic
with an isolated singularity at the origin.

In the light of substantial recent progress in the classification
of $2$--connected $7$--manifolds with $\rG_2$--structures \cites{Crowley2014,Crowley2015,Crowley2015b}, it is 
 a natural
task to sort such CY links $(K_f,\varphi)$. The important $\Z_{48}$--valued invariant $\nu(\varphi)$  introduced by Crowley
and Nordstr\"{o}m
 \cite{Crowley2015b} allows us to distinguish such pairs, up to diffeomorphism of $K_f $
and homotopy of $\varphi$, but its definition  is non-constructive and it
requires  an \emph{ad hoc} spin coboundary $8$--manifold $W$ such that $K_f =\partial W$. In Section \ref{sec: eta invariant}, we show that this coboundary can be essentially taken to be a typical Milnor fibre, and we find an explicit formula for $\nu(\varphi)$  in terms of topological data:

\begin{theorem}
\label{thm: Crowley-Nordstrom}
Let $K_f  \stackrel{S^1}{\longrightarrow} V\subset\P^4(w)$ be a weighted Calabi-Yau link of degree $d$ and weight  $w=(w_0,\dots,w_4)$; then  the Crowley-Nordstr\"om $\nu$ invariant of any $S^1$--invariant $\rG_2$--structure    $\varphi$
on $K_f $ is an odd integer given by $$
\nu(\varphi)=(\frac{d}{w_{0}}-1)\dots(\frac{d}{w_{4}}-1)-3(\mu_{+}-\mu_{-})+1
$$
where  $(\mu_{-},  \mu_{+})$ is the
 signature of the intersection form on $H^{4}(\widetilde{\cV},\R)$, for 
 $$
 \widetilde{\cV}={\{f(z)=1\}}\subset \mathbb{C}^5.$$
\end{theorem}
 Using a method by Steenbrink to calculate the signature, we obtain an
effective algorithm to determine $\nu(\varphi)$ for any CY link, with straightforward computational methods (see Appendix \ref{app: algorithm}). We observe that several values of $\nu$ are realised in this manner, and conjecture that indeed all possible 24 values can be realised by the `natural' cocalibrated $\rG_2$-structure (\ref{eq-GrayG2structure}) of a weighted CY link. In particular, for the homogeneous case we find: 
\begin{corollary}
\label{cor: nu for Fermat quintic}
The Crowley-Nordstr\"{o}m  $\nu$ invariant of the Fermat quintic link
with $\rG_2$--structure (\ref{eq-GrayG2structure}) is $\nu(\varphi)=5$.
\end{corollary}
To the best of our knowledge, this large class of $7$--manifolds with $\rG_2$--structure of the form $(K_f,\varphi)$ is the first instance besides the original reference \cite{Crowley2015b} in which the $\nu$  invariant has been computed explicitly.
\subsection{Gauge theory on contact Calabi-Yau manifolds}
In Section \ref{sec: gauge theory}, we turn   to the second axis of interest in $\rG_2$--geometry, as a model
for $7$--dimensional gauge theory. Since that concept appeared in the Physics literature \cite{Corrigan1983}, physicists pursue an analogous definition of Witten's  topological quantum field theory  \cite{Witten1988} on  spaces with $\rG_2$--metrics \cite{Acharya1997}. Moreover, it was noticed  in \cite{Harvey1999} that the superpotential for M--theory compactifications on $\rG_2$--manifolds `counts' associative  $3$--manifolds (i.e. submanifolds calibrated by $\varphi$)  in the same way as the prepotential of type II strings counts holomorphic curves in CY $3$--folds. Mathematicians, on the other hand, following the seminal viewpoint of \cite{Donaldson1998}, expect the theory to culminate in a topological count of instantons, yielding an invariant for $7$--manifolds with a $\rG_2$--structure, in the same vein as the Casson invariant for flat connections over $3$--manifolds \cite{Donaldson2002}. At the current stage, however, major compactification issues remain and a more thorough analytical understanding might have to be postponed in favour of  exploring a good number of examples \cites{SaEarp2009,Walpuski2013,Clarke2014,SaEarp2014,SaEarp2015a,SaEarp2015b}. 

 We propose a consistent formulation of
elementary Yang-Mills theory on $7$--dimensional cCY manifolds. In Section \ref{sec: YM and CS},
we define a connection $A$ on a complex vector bundle $E\to K$ to be a\emph{ $\rG_2$--instanton} if $F_A\wedge\psi=0$,
where $\psi$ is the  $\rG_2$--structure $4$--form (cf.  \cites{Donaldson1998,Tian2000}), which characterises $A$
at first as a critical point of the Chern-Simons functional. In Section \ref{sec: Sasakian vb} we endow
$E$ with a suitable  holomorphic Sasakian vector bundle structure, following
 the framework of Biswas and Schumacher \cite{Biswas2010}, to obtain a notion of  \emph{Chern connection},  compatible at once with the holomorphic structure
and some Hermitian bundle metric (Proposition \ref{prop: Chern connection}). In Section \ref{sec: G2 = tHYM}, we verify that the $\rG_2$--instanton condition is exactly
equivalent to a  natural transverse Hermitian Yang-Mills condition (Lemma
\ref{lem: SD lifts to G2-instanton} and Corollary \ref{cor: SD lifts to G2-instanton}),
somewhat similarly to the classical identification of selfdual and
HYM connections on compact K\"ahler surfaces \cite[\S2]{Donaldson1990}. Furthermore, going over into Section \ref{sec: top energy bounds},  the cCY    $\rG_2$--geometry allows us to move past some generally expected difficulties in the presence of torsion: since in these cases $\de\varphi=\omega\wedge \omega$,  intergable   $\rG_2$--instantons are  indeed \emph{Yang-Mills solutions}, i.e. $\de_A^*F_A=0$, and we define a secondary characteristic class leading to topological energy bounds, thus proving:   

\begin{theorem}
\label{thm: G2-inst are YM minima}
Let $\cE$ be a holomorphic Sasakian bundle over a $7$--dimensional closed contact Calabi-Yau manifold
$M$  endowed with its natural cocalibrated $\rm G_2$--structure (cf. Proposition \ref{prop: G2-structure on cCY}).
Then\begin{enumerate}
        \item 
         integrable $\rG_2$--instantons [cf. (\ref{eq: g2-instanton})] on $\cE$  
are critical points of the Yang-Mills functional $\cS_{\rm YM}$;
        \item
        if the absolute topological minimum of $\cS_{\rm YM}$  is attained among integrable connections, then the minima are exactly the $\rG_2$--instantons, i.e.,  the critical points of the Chern-Simons functional $\cS_{\rm CS}$.
\end{enumerate}
\end{theorem}

To conclude with an example,
in Section \ref{sec: example pullback} we focus on the simplest case in which $E$ is a pullback from the basic $CY^3$, and we derive the explicit local equations of $\rG_2$--instantons in that setting:       

\begin{theorem}
\label{thm: G2-intantons on K}
Suppose  $\pi: K\to V$ is a $7$--dimensional CY link, and let   $\cE:=\pi^*\cE_0\to K$ be
the pullback from a Hermitian  vector bundle $\cE_0\to V$. Then
\begin{enumerate}[(i)]
\item
        if an integrable connection  $\mathbf{A}=A+\sigma\theta$ on $E$ is  a  $\rG_2-$instanton,
then $A$ defines locally  a family $\left\{A_t\right\}_{t\in S^1}$ of Hermitian Yang-Mills connections on $\cE_0$, satisfying 
\begin{displaymath}
        \left(\frac{\partial A_t}{\partial t}-\de_{A_t}\sigma\right)\wedge\theta=0;
\end{displaymath} 
        \item
        if $\cE$ is indecomposable, there is a one-to-one correspondence between $S^1$--invariant
         $\rG_2$--instantons on $\cE$ and Hermitian Yang-Mills connections
         on $\cE_0$. 
\end{enumerate}
\end{theorem}

In particular, Theorem \ref{thm: G2-intantons on K} implies that  $S^1$--invariant
         $\rG_2$--instantons on $\cE$ are `counted' by the Donaldson-Thomas
         invariants of $\cE_0$, and this count should remain constant at least for
         any $S^1$--invariant deformations of the $\rG_2$--structure (\ref{eq-GrayG2structure}).
Finally, we underscore that the homogeneous case is offered as proof of principle, since our narrative seems to readily  extend to crepant resolutions of weighted projective Calabi-Yau $3$--orbifolds. 

Readers interested in a more detailed account of instanton theory on   $\rG_2$--manifolds are kindly
 referred to the introductory sections of \cites{SaEarp2015a,SaEarp2015b} and citations therein. A more thorough study of the moduli spaces of $\rG_2$--instantons, Sasakian HYM connections and also contact instantons on contact Calabi-Yau $7$--manifolds will appear shortly as part of Luis Portilla's PhD thesis \cite{Portilla2020}.

\bigskip         

\noindent\textbf{Ackowledgements:} We thank Gon\c calo Oliveira and Johannes
Nordstr\"om for contributions to the proof of Theorem \ref{thm: Crowley-Nordstrom}
and the Mathematical Sciences Research Institute in Berkeley, CA for hosting
that discussion. We also thank Thomas Walpuski for pointing out a mistake
in the preprint version. We are grateful to Sara D. Cardell for designing the figures. OC-A was supported by S\~ao Paulo
State Research Council (Fapesp) grant 2014/23594--6 and CONACYT 262121. LRD was supported by Fapesp grant 2014/13357-7 and HSE was supported by Fapesp grant 2014/24727-0 and by Brazilian National Research Council (CNPq) grant 312390/2014-9.  


\section{Geometric structures on links}
\label{sec-geometric-structures}

We address the possibilities of    $\rG_2$--geometry on Calabi-Yau links, starting from the motivational fact that a $7$--manifold  admits a $\rG_2$-structure if and only if it is orientable and spin  \cite{Gray1969},
as is the case of links weighted homogeneous hypersurface singularities in $\C^5$ \cite[Theorem $9.3.2$]{Boyer2008}. Such links have a very rich tautological geometry, including
a null Sasakian structure with a compatible non-degenerate  $3$--form which is `transversely' holomorphic, fitting in the category of contact Calabi-Yau manifolds proposed by Tomassini and Vezzoni \cite{Tomassini2008}. In this section we compile relevant definitions and known properties of weighted homogeneous links, and derive some straightforward consequences.      

\subsection{Hypersurface links of isolated singularities}

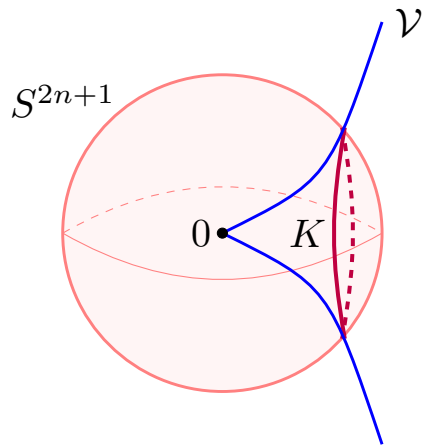
\begin{figure}[h]
\begin{tikzpicture}[scale=0.7]


\draw [line width=1.1,color=red!50,fill=red!4](0,0)circle(3);
\draw[dashed,color=red!50] (-3,0) to [bend left] (3,0);
\draw  [color=red!50](-3,0) to [bend right] (3,0);

\draw[line width=1.5,color=purple](2.3,2) to [bend right=10](2.3,-2);
K
\draw[line width=1.5, dashed, color=purple](2.25,2) to [bend left=10](2.25,-2);
K suspensiva

\draw [line width=1.1,blue] (0,0) .. controls (2,1) and (2,1) .. (3,4);
superior
\draw [line width=1.1,blue](0,0) .. controls (2,-1) and (2,-1) .. (3,-4);
inferior

\draw (0,0) node[draw, circle,scale=0.4,fill=black] {};
\draw(-0.4,0) node[scale=1.5] {$0$};

\draw (3.5,4) node[scale=1.5]{$\mathcal{V}$};
\draw (1.6,0) node [scale=1.5] {$K$};
\draw (-3,2.5) node [scale=1.5] {$S^{2n+1}$};
\end{tikzpicture}

\caption{Link $K$ of a hypersurface $\cV$ in $\C^{n+1}$}
\label{fig: hypersurface singularity link}
\end{figure}

We begin by reviewing more carefully Milnor's fibration theorem, following
the original reference \cite[\S 5-7]{Milnor1969}.
We denote by $\overline{B}_{\epsilon}$ the closed ball of radius $\epsilon$ centered at the origin of $\C^{n+1}$, by $S^{2n+1}_{\epsilon}=\partial B_{\epsilon}(0)$ the boundary of this ball, and $B_{\epsilon}$ for the corresponding open ball. Let $f:\C^{n+1}\rightarrow\C$ be a complex analytic map with $f(0)=0$ and denote $\cV:=f^{-1}(0)$ and $K_f :=\cV\cap S^{2n+1}_{\epsilon}$ (Figure
\ref{fig: hypersurface singularity link}).

\begin{theorem}
\label{thm-Milnorfibration}
Let $\epsilon > 0$ be sufficiently small; then the map $$\phi: S^{2n+1}_{\epsilon}-K\rightarrow S^{1}, \;\;\;\; \phi=\frac{f(x)}{\left|f(x)\right|},$$ is a locally trivial fibration, each fibre $F=\phi^{-1}(a)$ is smooth parallelisable and has the homotopy type of a finite CW-complex of dimension $n$. Furthermore, if $f$ has an isolated singularity at $0$, then each fibre $F$ has the homotopy type of a bouquet $S^{n}\vee \overset{\mu}{\dots} \vee S^{n}$ of spheres, and it is homotopy-equivalent to its closure $\overline{F}$ which is a compact manifold with boundary, with common boundary $\partial\overline{F}=K$. Likewise, $K_f $ is a smooth $(n-2)$-connected real manifold of dimension $2n-1$.
\end{theorem}
The number $\mu$ of spheres $S^{n}$ in the bouquet described in Theorem \ref{thm-Milnorfibration} is called the \emph{Milnor number} and it is an extremely  important topological invariant of the link.

\begin{theorem}\label{thm-Milnornumber}
The Milnor number $\mu$ has the following interpretations:
\begin{enumerate}[(i)]
        \item 
        $\mu$ is the complex dimension of the vector space obtained by taking the quotient of the local ring  $\cO_{0}(\C^{n+1})$  of holomorphic functions
at $0\in \C^{n+1}$ by the Jacobian ideal $J_{f}=\left(\partial f/ \partial z_{0}, \dots, \partial f/ \partial z_{n}\right)$ of $f$:
$$\mu=\mathrm{dim}_{\C}\frac{\cO_{0}(\C^{n+1})}{J_{f}},$$
        \item
        $\mu$ is the rank of the free Abelian middle homology group $H_{n}(F)$,
        \item
        $\mu$ is determined by the Euler characteristic of $F$: 
$$\chi(F)=1+(-1)^{n}\mu.$$
\end{enumerate}
\end{theorem}
 
The following result \cite[Theorem 5.11]{Milnor1969} gives a useful alternative description of the Milnor fibre:

\begin{theorem}[]\label{thm-Lefibration}
If a complex number $c\neq 0$ is sufficiently close to zero, then the complex hypersurface $f^{-1}(c)$ intersects the open ball $B_{\epsilon}$ along a smooth manifold which is diffeomorphic to the fibre $F$.
\end{theorem}

Now we focus on the particular case in which $f$ is a weighted homogeneous polynomial with an isolated singularity at $0\in \C^{n+1}$. This case is special because $\cV:=f^{-1}(0)$ intersects transversally every sphere $S^{2n+1}_{r}$ around the origin. Recall the definition of a weighted homogeneous polynomial.

\begin{definition}\label{def-weighted polynomial}
A polynomial $f(z_{0}, \dots, z_{n})$ is called a \emph{weighted homogeneous polynomial} of degree $d$ and weights $\left(w_{0}, \dots, w_{n}\right)$ if  for any $\lambda\in \C^{*}$
\begin{equation}
f(\lambda^{w_0}z_{0}, \dots,\lambda^{w_{n}} z_{n})=\lambda^{d}f(z_{0}, \dots, z_{n}).
\end{equation}
\end{definition}
\noindent NB.: a homogeneous polynomial of degree $d$ is weighted homogeneous of weights $\left(1, \dots, 1\right)$.

\begin{proposition}[\cite{Milnor1970}, Theorem $1$]
\label{prop: wh Milnor number}
Let $f(z_{0}, \dots, z_{n})$ be a weighted homogeneous polynomial of degree $d$ and weights $\left(w_{0}, \dots, w_{n}\right)$ having an isolated singularity at the origin. Then the cohomology   $H_{n}(F,\Z)$ is  free Abelian of rank $\mu=(\frac{d}{w_{0}}-1)\dots(\frac{d}{w_{n}}-1)$.
\end{proposition}

The Milnor fibration associated to a weighted homogeneous polynomial can appear under a different dressing, as the following lemma shows (\cite[Lemma
9.4]{Milnor1969}; see also \cite[Chapter 3, exercises 1.11 and 1.13]{Dimca1992}).

\begin{lemma}\label{lemm-affineMilnorfiber}
Let $f(z_{0}, \dots, z_{n})$ be a weighted homogeneous polynomial. Then the mapping $$f:\C^{n+1} - \cV\rightarrow \C^{*}$$ given by restriction of $f$ is a locally trivial fibration. Denote by $\psi$ the restriction of the above fibration over the unit circle $S^{1}$, then $\psi$ is fibre-diffeomorphic to the Milnor fibration $\phi$ of Theorem \ref{thm-Milnorfibration} associated to $f$. In particular the Milnor fibre is diffeomorphic to the non-singular affine hypersurface $\tilde\cV :=\{z\in\C^{n+1}| f(z)=1\}$.
\end{lemma}

Weighted homogeneous polynomials give rise in a natural way to links, fibering by circles over weighted projective hypersurfaces:
\begin{definition}
\label{def: wh link}
Let $f:\C^{n+1}\to\C$ be a $w$--weighted homogeneous polynomial with an isolated critical point at $0$, so that each sphere $S^{2n+1}=\partial B_\epsilon(0)$ intersects $\cV:=f^{-1}(0)\subset \C^{n+1}$ transversely.
Then $K_f := \cV\cap S^{2n+1}$ is called a \emph{weighted  link of degree
$\deg f$ and weight $w$}. 
\end{definition}
Given a weight vector $w=(w_0,\dots,w_n)$, denote by $\C^{*}(w)$ the weighted $\C^{*}$--action
on $\C^{n+1}$ given by
\begin{equation}\label{eq-weigted action}
(z_{0}, \dots, z_{n}) \rightarrow (\lambda^{w_0}z_{0}, \dots,\lambda^{w_{n}}
z_{n}).
\end{equation}
We have the commutative diagram 
\[      
\xymatrix@C1pc{
K \ar[r]^-{} \ar[d]_-{}& S^{9} \ar[d]^-{}\\
                V \ar[r]^-{}& \P^{4}(w),}
\]
where the horizontal arrows are Sasakian and K\"ahlerian embeddings, respectively,
and the vertical arrows are principal $S^{1}$--orbibundles and orbifold Riemannian
submersions. As a complex orbifold, the hypersurface $V\subset\P^4(w)$ can
be represented as the quotient $\left(\cV-{0}\right)/\C^{*}(w)$ where $\cV=f^{-1}(0)$.

\subsection{$\rG_2$-geometry}
\label{sec: G_2-geometry of links}
We now address the context of Theorem \ref{thm: cocalibrated G2-structure},
concerning the natural cocalibrated $\rG_2$--structure (\ref{eq-GrayG2structure}).
This section serves the double purpose of recalling notions of  $\rG_2$-geometry and setting the scene for the gauge theoretical investigation in Section \ref{sec: gauge theory}. 

Let $Y$ be an oriented smooth $7-$manifold. A $\rG_2-$\emph{structure} is a smooth tensor $\varphi \in \Omega ^{3}(Y)$ identified, at every  $p\in Y$, by some frame $f_{p}:T_{p}Y\rightarrow \mathbb{R}^{7}$, with the model
(sign conventions of \cite{Salamon1989})  
\begin{equation}
\label{eq: G2 3-form}
         \varphi_0 
        =e^{567}+\omega_1\wedge e^{5}
        +\omega_2\wedge e^{6}
        +\omega_3\wedge e^{7}
\end{equation}
in the sense that   $\varphi  _{p}=f_{p}^{\ast }\left(
\varphi _{0}\right)$,  where
\begin{displaymath}
        \omega_1= e^{12} - e^{34}, \quad 
        \omega_2= e^{13} - e^{42}, \quad 
        \text{and}\quad
        \omega_3= e^{14} - e^{23}
\end{displaymath}
are the canonical generators of selfdual $2$--forms in $\Lambda^2_+\left(\R^4\right)^*$.
The pointwise inner-product
\begin{equation*}        \label{eq: phi0 gives inner product}
        \left\langle u,v\right\rangle e^{1...7}
        :=\frac{1}{6}
        \left( u\lrcorner \varphi_{0}\right) 
        \wedge 
        \left( v\lrcorner \varphi _{0}\right) 
        \wedge 
        \varphi _{0}
\end{equation*}
determines a Riemannian metric $g_\varphi$ on $Y$,  under which $\ast_\varphi\varphi$ is given pointwise by
\begin{equation}
\label{eq: G2 4-form}
        \ast \varphi_0 
        =
        e^{1234}-\omega_1\wedge e^{67}
        -\omega_2\wedge e^{75}
        -\omega_3\wedge e^{56}
        .
\end{equation}
In the language of calibrated geometry \cite{Harvey1982}, a $7$--manifold with $\rG_2$--structure $(Y,\varphi)$ is said to be \emph{calibrated} if $\de\varphi =0$ and \emph{cocalibrated} if $\de\ast _{\varphi }\varphi =0$;
moreover it is common to omit $Y$ and refer simply to $\varphi$ in those terms. Cocalibrated $\rG_{2}$-structures appear in the Fern{\'a}ndez-Gray classification \cite{Fernandez1982} of $\rG_{2}$-structures by their intrinsic torsion. A  $\rG_2$--structure
$\varphi$ is both calibrated and cocalibrated if and only if  $\nabla^{g_\varphi}\varphi=0$,
in which case $\Hol(g_\varphi)\subseteq\rG_2$ 
and it is said to be \emph{torsion-free} \cite[Lemma 11.5]{Salamon1989}. 

Let us consider the following familiar example found, for example, in \cite[Proposition
11.1.2]{Joyce2000}:

\begin{example}
\label{ex: localmodel}
Let $\left(Z, \omega, \epsilon\right)$ be a Calabi-Yau 3-fold. Then the product manifold $Z\times S^{1}$ has a natural torsion-free $\rG_{2}$-structure defined by:
$$\varphi:= \de t\wedge \omega + \operatorname{Im}\epsilon,$$
where $t$ is the variable in $S^{1}$ and tensors are denoted identically
to their pullbacks under projection onto the $Z$ factor. The Hodge dual of $\varphi$ is 
$$\psi :=*\varphi=\frac{1}{2}\omega\wedge\omega + \de t\wedge\Re\epsilon$$
and the induced metric $g_{\varphi}=g_Z + \de t\otimes \de t$ is the Riemannian product metric on $Z\times S^{1}$, with holonomy $\operatorname{Hol}(g_\varphi)=\operatorname{SU}(3)\subsetneq\rG_{2}$.

\end{example}

In the case of a CY link, we only deviate from the product model
of Example \ref{ex: localmodel} in the sense that  $K_f $ is necessarily \emph{nontrivial} as a circle bundle over the $CY^3$ base $V$, since $\pi_1(K)=\{1\}$ by Theorem \ref{thm-Milnorfibration},
so it is fair to ask whether $K_f $ also inherits a `globally twisted'   $\rG_2$--structure
from the Calabi-Yau structure of $V$.  

\begin{remark}
\label{rem: decomposition}
If a Lie group $G$ induces a $G$-structure on a manifold $M$, then every bundle of tensors splits into summands corresponding to  irreducible representations of
$G$. The link $K_f $ carries a $G_2$-structure so, in particular, $2$--forms
split as
$$
\Omega^{2}(K)=\Omega^2_7(K)\oplus\Omega^2_{14}(K),
$$
where $\Omega^2_7(K)$ and $\Omega^2_{14}(K)$ are vector subbundles of $\Omega^{2}(K)$ with fibres isomorphic to the irreducible \textbf{7} and \textbf{14}  representations of $\rG_2$, respectively. It is a well-known fact about manifolds with a
$\rG_2$--structure  \cite{Bryant1987,SaEarp2015a}
that $\left( \Omega^{2}\right)_{^{\;7}_{14} }$ is respectively the $_{+1}^{-2}-$eigenspace of the  $%
\rG_{2}-$equivariant linear map%
\begin{eqnarray*}
        T_{\varphi} \;: \; \Omega ^{2} &\rightarrow& \Omega ^{2} \\
        \eta &\mapsto &T_{\varphi}\eta :=\ast \left( \eta \wedge \varphi
        \right) .
\end{eqnarray*}
\end{remark}

\subsection{Links as Sasakian 7--manifolds}

A  \emph{contact manifold}  $(M,\theta)$ is given by a smooth  $(2n+1)$--manifold $M$ and a \emph{contact structure}  $\theta\in\Omega^1(M)$ such that $\theta \wedge \left( \de \theta \right)^{n}\neq 0$, everywhere on $M$. On a contact manifold
there exists a unique \emph{Reeb vector field} $\xi\in\Gamma(TM)$, such that $\xi \lrcorner \theta=1$ and $\xi\lrcorner \de\theta=0$. The Reeb vector field is nowhere-vanishing, so it uniquely determines a $1$--dimensional foliation $N_{\xi}$  called the \emph{characteristic foliation}. It is customary
to think of contact manifolds as odd-dimensional analogues of symplectic
manifolds, with the $2$--form $\de \theta$ being `transversely symplectic'
with respect to the characteristic foliation. From that perspective, Sasakian
geometry encodes the notion of `transversely K\"ahler structure': 

\begin{definition}
\label{def: Sasakian manifold}
A \emph{Sasakian structure} on $M$ is a quadruple $\left(M, \theta, g, \Phi \right)$ such that $\left(M, g\right)$ is a Riemannian manifold, $(M,\theta)$ is a contact manifold with Reeb vector field $\xi$,  $\Phi$ is a global section
of $\End(TM)$, and the following relations hold:
\small
\begin{align*}
g\left(\xi,\xi\right)=1, \quad \Phi\circ \Phi=-\operatorname{Id}_{TM}+\theta \otimes \xi, \quad g\left(\Phi X,\Phi Y\right)=g\left(X,Y\right)-\theta\left(X\right)\theta\left(Y\right),\\
\nabla^{g}_{X}\xi=-\Phi X, \quad \left(\nabla^{g}_{X}\Phi\right)\left(Y\right)=g\left(X,Y\right)\xi-\theta\left(Y\right)X, \quad \qquad \qquad
\end{align*}
\normalsize
where $X, Y$ are vector fields on $M$ and $\nabla^{g}$ is the Levi-Civita connection corresponding to $g$.
If $\left(M, \theta, g, \Phi \right)$ satisfies these conditions we say  $M$ is a \emph{Sasakian manifold}.
\end{definition}

 If the orbits of  $\xi$ are all closed, hence circles, then $\xi$ integrates
to an isometric $\rU(1)$ action on $M$, in particular this action is locally free. If the action is in fact free then the Sasakian structure is said to be \emph{regular}, otherwise, it is said to be \emph{quasi-regular}.
The leaf space $\mathcal{Z}: = M/N_{\xi} = M/\rU(1)$ has the structure of a manifold or orbifold, in the regular or quasi-regular case respectively.

The sphere $S^{2n+1}$ has a natural contact structure given by the Hopf fibration $S^{2n+1}\stackrel{\pi}{\rightarrow} \P^n$: 
$$
\theta_{c}=-\frac{\bi}{2}\sum^{n}_{j=0}\left(z_{j}d\bar{z}_{j}-\bar{z}_{j}d z_{j}\right)
=
\sum^{n}_{j=0}\left(y_{j}dx_{j}-x_{j}dy_{j}\right)
$$ 
(in the real coordinates $z_{j}=x_{j}+\bi y_{j}$).
It carries moreover \cite{Sasaki1962} a regular Sasakian structure $\left(S^{2n+1}, \theta_{c}, g_{c}, \Phi_{c} \right)$
in which the Reeb vector field is  
\small
$$\xi_{c}=\sum_{i=0}^{n}\left(y_{i}\frac{\partial}{\partial x_{i}}-x_{i}\frac{\partial}{\partial y_{i}}\right)=-\bi\sum_{j=0}^{n}\left(z_{j}\frac{\partial}{\partial z_{j}}-\bar{z}_{j}\frac{\partial}{\partial \bar{z}_{j}}\right),$$
\normalsize
the metric $g_{c}$ is given by the  inclusion $S^{2n+1}\subset \R^{2n+2}$, and

$$
\Phi_{c}=\sum_{i,j}\left\{
\left[\left(x_{i}x_{j}-\delta_{ij}\right)\partial_{x_i}
+
\left(x_{j}y_{i}\right)\partial_{y_i}
\right]
\otimes dy_{j}
-
\left[\left(y_{i}y_{j}-\delta_{ij}\right)\partial_{y_i}
+
x_{i}y_{j}\partial_{x_i}
\right]
\otimes dx_{j}
\right\}.
$$

The links of isolated hypersurface singularities admit Sasakian structures in a natural way.

\begin{proposition}[{Proposition 9.2.2 \cite{Boyer2008}}]
Let $K_{f}$ be the link of a hypersurface singularity. Then the Sasakian structure $\mathcal{S}_{c}:=\left(\theta_{c}, g_{c}, \Phi_{c} \right)$ on $S^{2n+1}$ defined above induces by restriction a Sasakian structure, also denoted by $\mathcal{S}_{c}$, on the link $K_{f}$.
\end{proposition}

\subsection{Contact Calabi-Yau structures on links}
\label{sec: contact Calabi Yau}
Contact Calabi-Yau manifolds were introduced by Tomassini and Vezzoni in \cite{Tomassini2008} and thoroughly studied by Habib and Vezzoni in \cite{Habib2015},
as a development of Reinhart's general theory of Riemannian foliations \cite{Reinhart1959}. This concept
describes Sasakian manifolds endowed with a closed basic complex volume form, which is `transversally holomorphic' in a certain sense (see Definition \ref{def: contact CY mfd}).
Most importantly for us, it allows for a vast generalisation of the $\rG_2$--geometry on homogeneous links discussed in Section \ref{sec: G_2-geometry of links}.

Let $(M, \theta)$ be a contact manifold with contact $1$-form $\theta$ and denote $B:= \operatorname{ker}\theta$ its contact distribution of rank $2n$, i.e., $TM=B\oplus N_{\xi}$. Let $X$ denote an arbitrary vector field tangent to the characteristic foliation $N_{\xi}$. A differential  form $\beta\in\Omega^k(M)$ is said to be \emph{transversal} if $X \lrcorner \beta =0$ and $\mathcal{L}_{X} \beta=0$ for every such  $X$. 

If $(M, \theta, \Phi)$ is a Sasakian manifold, and $x\in M$, it follows from Definition \ref{def: Sasakian manifold} that $\left(\Phi|_{B_{x}}\right)^{2}=-\operatorname{Id}_{B_{x}}$. Then we can decompose the complexification $B_{x}\otimes_{\mathbb{R}}\mathbb{C}$ into the eigenspaces of the complexifed automorphism  $\Phi|_{B_{x}}\otimes_{\mathbb{R}}\mathbb{C}$:
$$B_{x}\otimes_{\mathbb{R}}\mathbb{C}=B_{x}^{1,0}\oplus B_{x}^{0,1},$$
where $B_{x}^{1,0}$ and $B_{x}^{0,1}$ correspond to the eigenvalues $\bi:=\sqrt{-1} $ and $-\bi $ respectively. This induces a splitting of the exterior differential algebra over $B_\C :=B\otimes_{\mathbb{R}}\mathbb{C}$:
\begin{equation}
\label{eq: B^(p,q)}
\Omega^{k}\left( B_\C \right)=\bigoplus_{p+q=k}\Omega^{p,q}(M),
\end{equation}
where $\Omega^{p,q}(M):=\Gamma\left(\Lambda^{p}(B^{1,0})^*\otimes (\Lambda^{q}B^{0,1})^*\right)$ and $p,q \geq 0$. Then we have an obvious  decomposition of exterior forms
on $M$ given by
$$
\Omega^j(M)=\bigoplus_{p+q=j}\Omega^{p,q}(M)
\oplus
\bigoplus_{p+q=j-1}\Omega^{p,q}(M)
\wedge
\theta.
$$
If $\beta\in\Omega^k(M)$ is a transversal differential form, we will say that \emph{$\beta$ is of type $(p,q)$} if it belongs to $\Omega^{p,q}(M)$. The following lemma \cite[Corollary 3.1]{Biswas2010} will be crucial for our applications in gauge theory, so we sketch for convenience:

\begin{lemma}
\label{lemma: dtheta is (1,1)}
Let $(M, \theta, \Phi)$ be a Sasakian manifold. Then $d\theta \in \Omega^{1,1}(M)$.
\begin{proof}
That $d\theta$ is transversal is clear from Definition \ref{def: Sasakian manifold} . It is easy to prove that $d\theta(X,Y)=-g(\Phi X, Y)$ for all $X, Y\in B$, then $d\theta$ is of type $(1,1)$. 
\end{proof}
\end{lemma}

\begin{definition}
\label{def: contact CY mfd}
A \emph{contact Calabi-Yau manifold (cCY)} is a quadruple $(M, \theta, \Phi, \epsilon)$ such that:
\begin{itemize}
\item $(M, \theta, \Phi)$ is a $2n + 1$-dimensional Sasakian manifold;
\item $\epsilon$ is a nowhere vanishing transversal form on $B=\operatorname{ker}(\theta)$ of type $(n,0)$:
\begin{gather*}
\epsilon \wedge \bar{\epsilon}=c_{n}\omega^{n}, \quad d\epsilon=0, 
\end{gather*}
where $c_{n}=(-1)^{n(n+1)/2}\bi^{n}$ and $\omega:=d\theta$.
We denote accordingly
$$
\Re\epsilon:=\frac{\epsilon+\bar{\epsilon}}{2}
\qandq
\Im \epsilon:=\frac{\epsilon-\bar{\epsilon}}{2\bi}.
$$
\end{itemize}
\end{definition}

Our interest in cCY structures for  $\rm G_2$-geometry derives from the following fundamental result:

\begin{proposition}[\cite{Habib2015}, subsection $6.2.1$]
\label{prop: G2-structure on cCY}
Let $(M, \theta, \Phi, \epsilon)$ be $7$-dimensional contact Calabi-Yau manifold. Then $M$ carries a cocalibrated $\rm G_2$--structure defined by
\begin{equation}
\label{eq: cCY G2-structure}
        \varphi:= \theta\wedge \omega + \Im\epsilon
\end{equation} 
with torsion $d\varphi=\omega\wedge\omega$ (cf. Definition \ref{def: contact CY mfd}) and corresponding dual $4$-form
$$
\psi=\ast\varphi=\frac{1}{2}\omega\wedge\omega+\theta\wedge\Re\epsilon 
$$
\end{proposition}

The existence of cCY structures on links is equivalent to a simple numerical criterion on the weighted homogeneous data, which we adopt as a definition:

\begin{definition}
\label{def: CY link}
A weighted  link $K_f $ (cf. Definition \ref{def: wh link}) of degree $d$
and weight $w=\left(w_{0}, \dots, w_{n}\right)$ is said to be a \emph{Calabi-Yau
 (CY) link} if 
$$
d=\sum_{i=0}^{n}w_{i}.
$$
\end{definition}

The condition $d-\sum_{i=0}^{n}w_{i}=0$ means precisely that the Sasakian structure $(K, \theta_{c}, \Phi_{c})$ on $K_f $ induced from the canonical Sasakian structure of the sphere $S^{2n+1}$ is null Sasakian, i.e., the (basic) first Chern class of $(K, \theta_{c}, \Phi_{c})$ vanishes. Recall also this vanishing is exactly
the requirement for the weighted projective $V$ to be a Calabi-Yau orbifold \cite{Candelas1990}, thus CY links are nontrivial circle fibrations over Calabi-Yau $3$--orbifolds.
Furthermore, the Reeb vector field the unit tangent to the $S^1(w)$-action and the $3$-form $\varepsilon$ is transversal, so the  $\rm G_2$--structure
(\ref{eq: cCY G2-structure}) is $S^1$--invariant.
 In the terms of Definition \ref{def: CY link}, Habib and Vezzoni's existence result can be restated as:

\begin{proposition}[\cite{Habib2015}, Proposition $6.7$]
\label{prop: cCY link}
Every Calabi-Yau link admits a $S^1$--invariant contact Calabi-Yau structure.
\end{proposition}

The proof of Proposition \ref{prop: cCY link} relies on a Sasakian version of the El Kacimi theorem to prove that any null Sasakian structure on a compact simply-connected manifold can be 
deformed into a contact Calabi-Yau one.
Combining the previous two propositions:

\begin{corollary}[\cite{Habib2015}, Corollary $6.8$]
\label{cor: cocalibrated G2-str on K}
Every Calabi-Yau link has a cocalibrated $S^1$--invariant $\rm G_2$--structure of the form (\ref{eq: cCY G2-structure}).
\end{corollary}

\section{The  $\nu$ invariant of Calabi-Yau links}

\label{sec: eta invariant}

For an arbitrary closed $7$-manifold with  $\rm G_2$-structure $(Y^7,\varphi)$, Crowley and Nordstr\"om define a pair of homotopy invariants $(\nu(\varphi),\xi(\varphi))$
which completely classifies the data, up to diffeomorphism and homotopy, if
$Y$ is $2$--connected \cite[Theorem 1.17]{Crowley2015b}. Subsequently this
has been refined as an analytic invariant of manifolds with $\rG_2$--metrics
\cite{Crowley2015}, and similar ideas also intervene in the authors' topological
classification
of spin $2$--connected $7$--manifolds  \cite{Crowley2014}.

We will be interested in the first invariant $\nu(\varphi)$, which is a 
 $\mathbb{Z}_{48}$--valued combination of topological data from a compact
coboundary $8$--manifold with a $\rm{Spin}(7)$--structure $(W^8,\Psi)$ filling $(Y,\varphi)$,
in the sense that $Y=\partial W$ and $\Psi|_{Y}=\varphi$:
\begin{equation}
\label{eq: nu(phi)}
        \nu(\varphi):=\chi(W)-3\sigma(W) \quad\rm{mod}\; 48 
\end{equation}
($\chi$ and $\sigma$ denote the real Euler characteristic and
the signature, respectively.)
This quantity is preserved under diffeomorphisms of $Y$ and homotopies of the $\rG_2$--structure
$\varphi$ \cite[Theorem 1.3]{Crowley2015b}. Moreover, $\nu(\varphi)$ is independent of the particular choice of coboundary $W$ \cite[Corollary 3.2]{Crowley2015b},
thus it is interpreted as an ``$\hat A$-defect" from certain integral  characteristic
classes of principal $\rm{Spin}(8)$--bundles evaluated on $TW$ and the half-spinor
bundles $S^\pm W$. 

A central aspect is the fact that such a filling
$W$ always exists \cite[Lemma 3.4 (ii)]{Crowley2015b}. 
The argument 
relies on the fact that the bordism group $\Omega^{\rm{spin}}_7$ is trivial,
hence there always exists \emph{some}  (connected) coboundary $(W,\Psi)$ inducing a reference $\rm G_2$--structure on $Y$, but it is totally non-constructive. For example,  the authors must resort to an elaborate construction
of an explicit  coboundary $W$ to calculate $\nu=24$ [Theorem 1.7] for the important class of manifolds with
holonomy $\rm G_2$ obtained as \emph{twisted connected sums} \cite{Corti2015}.
This allows one to distinguish, for instance, whether a given $\rm G_2$--structure
is not a gluing of asymptotically cylindrical Calabi-Yau $3$-folds \cite{Corti2013}.
\subsection{Construction of a spin coboundary}
In order to calculate the $\nu$ invariant for our $\rG_2$--structure (\ref{eq: cCY G2-structure})
on a  link $K_f $, we must therefore find an ad hoc compact $\rm{Spin}(7)$--coboundary
$(W,\Psi)$ such that:
$$
K=\partial W \qandq \Psi|_{K}=\varphi.$$

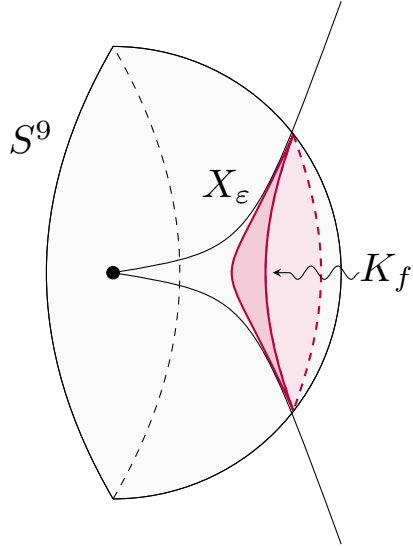
\begin{figure}[h]
\begin{center}
\begin{tikzpicture}[scale=0.6]
\draw[fill=gray!4] ([shift=(-90:5)]0,0) arc (-90:90:5)--(0,5) to [bend right]
(0,-5);

\draw[dashed]  (0,-5) to [bend right] (0,5);

\draw[thick, purple, dashed,fill=purple!10](3.95,3.05) to [bend left=20](3.95,-3.05) to [bend left=20](3.95,3.05);
\draw[thick, purple,fill=purple!20](3.95,3.05) to [bend right=20](3.95,-3.05)..controls
(2.15,1) and  (2.15,-1) ..(3.95,3.05);

\draw  (0,0) .. controls (3,0.5) and (3,0.5) .. (5,6);

\draw (0,0) .. controls (3,-0.5) and (3,-0.5) .. (5,-6);

\draw (0,0) node[draw, circle,scale=0.5,fill=black] {};

\draw ([shift=(-90:5)]0,0) arc (-90:90:5)--(0,5) to [bend right] (0,-5);
de la esfera y ecuador

\draw  (-1.8,3) node [scale=1.5]{$S^9$};
\draw  (2.5,1.9) node [scale=1.5]{$X_{\epsilon}$};
\draw  (6,0) node [scale=1.5]{$K_f$};

\draw[-stealth,decorate,decoration={snake,post=lineto, post length=2mm}] (5.4,0)--(3.5,0);
\end{tikzpicture}
\end{center}
\caption{Smoothing $X_\epsilon$ of $\cV$ inside $S^9$, with boundary $K_f$.}
\label{fig: smoothing}
\end{figure}

Let $K_f $ be the weighted  link (cf. Definition \ref{def: wh link}) of degree $d$ and weight $w=\left(w_{0}, \dots, w_{4}\right)$. 
The ambient  $4$--form
$$
\Psi:=\frac{1}{f}\sum_{i=0}^{4} z_i dz_0\wedge\overset{\hat i}{\dots}\wedge dz_4
\in \Lambda^{4,0}(\C^5)
$$
is $S^{1}\subset \C^{*}$--invariant under the action (\ref{eq-weigted action}) if, and only if, $d-\sum_{i=0}^{4}w_{i}=0$, i.e., exactly when the link $K_f $ is Calabi-Yau (Definition \ref{def: CY link}). 
Let $\feps$ be a smoothing of $f$, e.g. $\feps:=f-\epsilon$ and 
$$
X_\epsilon:=\feps^{-1}(0)\cap \overline{B}^{10}\subset \C^5
$$ 
the $8$--manifold inside the (compact component of the complement of the)
sphere $S^9$ with boundary $K_f =\partial\xeps$ (Figure \ref{fig: smoothing}). The restriction $\Psi|_{\xeps}$
induces an $SU(4)$--structure, hence a $\Spin(7)$--structure on $\xeps$,
which is $S^1$--invariant by construction. 

Restricting to the boundary we get an  $S^1$--invariant $\rG_2$--structure
on $K_f $, which corresponds exactly to an $SU(3)$--structure $\varphi'$ on
the basic Calabi-Yau $V$. Now, all   $SU(3)$--structures on a $6$--manifold
are homotopic, as sections of a bundle of rank $8$, so $\varphi'$ is homotopic
to our $\varphi$ and we can take $W=\xeps$: 
\begin{equation}
\label{eq: nu from smoothing X}
        \nu(\varphi)=\nu(\varphi')=\chi(\xeps)-3\sigma(\xeps) \mod\;48.
\end{equation}
Now all we need is to calculate the topology of the smoothing of an affine
hypersurface.  

\begin{proposition}
\label{prop: X = F = V}
Let $f:\C^{n+1}\to\C$ be weighted homogeneous polynomial with typical Milnor fibre  $F$ (cf. Theorem \ref{thm-Milnorfibration}), and consider the model affine variety
$$
\widetilde{\cV}:=\left\{ z\in\C^{n+1} : f(z)=1 \right\}.
$$
Given $\varepsilon>0$ sufficiently small, the smoothing $\feps:=f-\epsilon$, defines in $\C^{n+1}$ a compact manifold with boundary $X_\epsilon:=\feps^{-1}(0)\cap \overline{B}^{2n+2}$, and the following  are diffeomorphic:
\begin{gather}
\xeps\simeq \overline{F} \qquad \xeps \backslash \partial \xeps \simeq F \simeq \widetilde\cV.
\end{gather}
\begin{proof}
Taking $\epsilon=c=1$ in Theorem
\ref{thm-Lefibration}, we identify diffeomorphically  the smoothing $X_{\epsilon}$
with the closure of the Milnor fibre $\overline{F}$. Then the second identification is immediate from Lemma \ref{lemm-affineMilnorfiber}.
\end{proof}
\end{proposition}

\subsection{Explicit formula for $\nu$ on Calabi-Yau links}
\label{sec: Formula for nu}
In view of Proposition \ref{prop: X = F = V}, we will obtain the $\nu$ invariant from the following:
\begin{equation}
\label{eq: sigma and chi of X }
         \chi(\xeps)
        =\chi(\overline{F})
        \qandq
        \sigma(\xeps)=\sigma(\widetilde{\cV}).
\end{equation}

We begin with
Steenbrink's method  \cite{Steenbrink1977} for the signature of $\widetilde{\cV}$. 
 Let $\{z^{\alpha}: \alpha=(\alpha_{0}, \dots, \alpha_{n})\in I \subset \N^{n+1}\}$
be a set of monomials in $\C[z_{0}, \dots, z_{n}]$ representing a basis over
$\C$ for $\tfrac{\C[[z_{0}, \dots, z_{n}]]}{\left(\partial
f/\partial z_{0}, \dots, \partial f/\partial z_{n}\right)}$ (cf. (i) of Theorem \ref{thm-Milnornumber}). For each $\alpha\in
I$ define 
\begin{equation}\label{eq-array Steenbrink}
l(\alpha):=\sum^{n}_{i=0}\left(\alpha_{i}+1\right)\frac{w_{i}}{d}.
\end{equation}
Assume that $n$ is even (in our case, indeed $n=4$), and denote by $(\mu_{-}, \mu_{0}, \mu_{+})$ the
signature of the intersection form on $H^{n}(\widetilde{\cV},\R)$ i.e., $\mu_{-}$, $\mu_{0}$
and $\mu_{+}$ denote the numbers of negative, zero and positive entries,
respectively, on the diagonal of the intersection matrix. 
Then
$$
\sigma(\widetilde{\cV})=\mu_{+}-\mu_{-}.
$$
Note that the sum $\mu_{+}+\mu_{-}+\mu_{0}$ equals the Milnor number $\mu$ by (ii) of Theorem \ref{thm-Milnornumber}. On the other hand, by (iii) of Theorem \ref{thm-Milnornumber}, the Euler characteristic of the Milnor fiber is determined by  the Milnor number, which is given by Proposition \ref{prop: wh Milnor
number} for weighted homogenous links.  By Theorem \ref{thm-Milnorfibration},  $F$ is homotopy-equivalent
to $\overline{F}$, so for $n=4$:   
$$\chi(\overline{F})=\chi(F)=1+(\frac{d}{w_{0}}-1)\dots(\frac{d}{w_{4}}-1).$$
Finally, replacing (\ref{eq: sigma and chi of X }) in (\ref{eq: nu from smoothing X}), we establish the formula of Theorem \ref{thm: Crowley-Nordstrom}: 
\begin{equation}
\label{eq: nu(phi) = formula}
        \nu(\varphi)=(\frac{d}{w_{0}}-1)\dots(\frac{d}{w_{4}}-1)-3(\mu_{+}-\mu_{-}) +1.
\end{equation}

Steenbrink proved
\cite[Theorem~2]{Steenbrink1977} that  the signature $(\mu_{-}, \mu_{0}, \mu_{+})$ can be
computed as follows:
\begin{eqnarray*}
        \mu_{+}&=&\left\vert\{\beta \in I: l(\beta)\notin\Z, \left\lfloor
l(\beta)\right\rfloor \in 2\Z  \}\right\vert,\\
        \mu_{-}&=&\left\vert\{\beta \in I: l(\beta)\notin\Z, \left\lfloor
l(\beta)\right\rfloor \notin 2\Z \}\right\vert,\\
        \mu_{0}&=&\left\vert\{\beta \in I: l(\beta) \in\Z \}\right\vert,
\end{eqnarray*}
where $\left\lfloor x\right\rfloor$ denotes the integer part of $x\in\mathbb{Q}$, hence the above process can be easily implemented. We offer the code for a working algorithm
in a combination of \textsc{Singular} and \textsc{Mathematica}, but surely
 readers will be able to formulate leaner alternatives.
We display in Table \ref{tab: values of nu for various CY} the invariants given by (\ref{eq: nu(phi) = formula}) for some examples from Candelas' list of weighted Calabi-Yau threefolds. 

\begin{table}[h!]
  \centering
  \caption{The $\nu$ invariant for certain Calabi-Yau links}
  \label{tab: values of nu for various CY}
  \begin{tabular}{l|l|l|l}
        \toprule
        degree & weights & polynomial & $\nu$ \\
        \midrule                                
         75 & (10,12,13,15,25) & $z_{0}^{5}z_{4}+z_{1}^{5}z_{3}+z_{2}^{5}z_{0}+z_{3}^{5}+z_{4}^{3}$  & 1\\
                                \hline
         135 & (1,18,32,39,45) & $z_{0}^{135}+z_{1}^{5}z_{4}+z_{2}^{3}z_{3}+z_{3}^{3}z_{1}+z_{4}^{3}$ & 3\\
                                \hline
         36 & (18, 12, 4, 1, 1) &  $z_{0}^{2}+z_{1}^{3}+z_{2}^{9}+z_{3}^{36}+z_{4}^{36}$   & 5\\
         \hline
         81 &  (3,7,18,26,27)  &  $z_{0}^{27}+z_{1}^{9}z_{2}+z_{2}^{3}z_{4}+z_{3}^{3}z_{0}+z_{4}$  &  7 \\
        \hline
        45  &  (3,5,8,14,15)  &  $z_{0}^{15}+z_{1}^{9}+z_{2}^{5}z_{1}+z_{3}^{3}z_{0}+z_{4}^{3}$  &  9\\
                                \hline
                                45  & (4,7,9,10,15)  & $z_{0}^{9}z_{2}+z_{2}^{5}+z_{1}^{5}z_{3}+z_{3}^{3}z_{4}+z_{4}^{3}$ & 11\\
        \hline
        75 & (5,8,12,15,35)   & $z_{0}^{15}+z_{1}^{5}z_{4}+z_{2}^{5}z_{3}+z_{3}^{5}+z_{4}^{2}z_{0}$   &13\\
                                \hline
       180 & (90, 60, 20, 9, 1) &  $z_{0}^{2}+z_{1}^{3}+z_{2}^{9}+z_{3}^{20}+z_{4}^{180}$  & 15\\
                        \hline
        45 & (15, 15, 5, 9, 1) &  $z_{0}^{3}+z_{1}^{3}+z_{2}^{9}+z_{3}^{5}+z_{4}^{45}$   &   17 \\                     
        \hline
       16 & (4,8,2,1,1)   & $z_{0}^{2}z_{1}+z_{1}^{2}+z_{2}^{4}z_{1}+z_{3}^{16}+z_{4}^{16}+z_{2}^{8}$ & 19\\
        \hline
        81 &  (2,9,19,24,27) &  $z_{0}^{27}z_{4}+z_{2}^{3}z_{3}+z_{3}^{3}z_{1}+z_{1}^{9}+z_{4}^{3}$  &  21\\
                                \hline
        24 &  (12, 8, 2, 1, 1) & $z_{0}^{2}+z_{1}^{3}+z_{2}^{12}+z_{3}^{24}+z_{4}^{24}$  &  23 \\
                                \hline
        1806 & (42, 258, 903, 602, 1) &  $z_{0}^{43}+z_{1}^{7}+z_{2}^{2}+z_{3}^{3}+z_{4}^{1806}$   & 25 \\
                                \hline
                                51 & (2,6,9,17,17)  & $z_{0}^{17}z_{4}+z_{1}^{7}z_{2}+z_{2}^{5}z_{1}+z_{3}^{3}+z_{4}^{3}$  &  29\\
                                \hline
                                93 & (3,8,21,30,31) & $z_{0}^{31}+z_{1}^{9}z_{2}+z_{2}^{3}z_{3}+z_{3}^{3}z_{0}+z_{4}^{3}$  &  31\\
                                \hline
                                63 & (3,4,14,21,21) & $z_{0}^{21}+z_{1}^{15}z_{0}+z_{2}^{3}z_{3}+z_{3}^{3}+z_{4}^{3}$  & 33\\
                                \hline
                                103 & (1,16,23,29,34) & $z_{0}^{103}+z_{1}^{5}z_{2}+z_{2}^{3}z_{4}+z_{3}^{3}z_{1}+z_{4}^{3}z_{0}$ & 37\\
                                \hline
                                135 & (5,6,14,45,65)  & $z_{0}^{27}+z_{4}^{2}z_{0}+z_{1}^{15}z_{3}+z_{3}^{3}+z_{2}^{5}z_{4}$ &  39\\
                                                                                                                                \hline
                                                                                                                                60 & (30, 20, 5, 4, 1) & $z_{0}^{2}+z_{1}^{3}+z_{2}^{12}+z_{3}^{15}+z_{4}^{60}$   &  41 \\
                                \hline
                                55 & (4,4,11,17,19)  & $z_{0}^{11}z_{2}+z_{1}^{9}z_{4}+z_{2}^{5}+z_{3}^{3}z_{1}+z_{4}^{2}z_{3}$  & 43\\
                                \hline
                                135 & (1,21,30,38,45) & $z_{0}^{135}+z_{1}^{5}z_{2}+z_{2}^{3}z_{4}+z_{3}^{3}z_{1}+z_{4}^{3}$  &  45\\
                                                                                                                                \hline
                              45 & (5, 5, 9, 11, 12)  & $z_{0}^{9}+z_{1}^{8}z_{0}+z_{2}^{5}+z_{4}^{3}z_{2}+z_{3}^{3}z_{4}$  & 47 \\
                             \bottomrule
  \end{tabular}
\end{table} 
 
Inspection of a few examples suggests a  parity constraint for the $\nu$ invariant, and this is indeed the case:

\begin{proposition}
The Crowley-Nordstr\"om  $\nu$ invariant of a weighted  link is odd in $\Z_{48}$.
\begin{proof}
We know from \cite[Theorem 1.3]{Crowley2015b} that $\nu(\varphi)\equiv \chi_{\Q}(K)\;\rm{mod}\;2$, where 
$$
\chi_{\Q}(K):=\sum_{i=0}^{n-1}b_{i}(K)
$$ is the rational semi-characteristic of $K_f $. On the other hand, $b_1=b_2=0$ because $K_f $ is $2$--connected
(cf. Theorem \ref{thm-Milnorfibration}),
and we know from  \cite[Theorem $9.3.2$]{Boyer2008}
that the Betti number $b_{n-1}$ is even, if $n$ is even. Therefore   $b_{3}$ is even when $n=4$, thus $\chi_{\Q}(K)$
is odd.
\end{proof} 
\end{proposition}
Together with formula (\ref{eq: nu(phi) = formula}), this completes the proof of Theorem \ref{thm: Crowley-Nordstrom}.

\begin{example}
\label{exa: Fermat quintic}
Let us calculate the $\nu$ invariant for our $\rG_2$--structure (\ref{eq-GrayG2structure})
on  the Fermat quintic 
$$
f(z)=z_0^5+z_1^5+z_2^5+z_3^5+z_4^5.
$$
In this case the Milnor algebra is just $\tfrac{\C[z_{0}, \dots, z_{4}]}{\left(z_{0}^{4}, \dots, z_{4}^{4}\right)}$, the Milnor number is $\mu=1024$ and we can take, as a basis of the Milnor algebra, all monomials of the form $z_{0}^{\alpha_0} \dots z_{4}^{\alpha_4}$, with $0\leq \alpha_{i} \leq 3,\;\forall i$. A simple computation gives 
$$
(\mu_{-}, \mu_{0}, \mu_{+})=(240, 204, 580),
$$ 
therefore $\sigma(\widetilde{\cV})=340.$  On the other hand, Theorem \ref{thm-Milnornumber} gives $\chi(\xeps)
=1+\mu=1025.$
Hence $\nu(\varphi)=5$ as claimed in Corollary \ref{cor: nu for Fermat quintic}. 
\end{example}

\section{Gauge theory on contact Calabi-Yau $7$--manifolds}
\label{sec: gauge theory}

Let $(M,\varphi)$ be a closed contact Calabi-Yau manifold and consider a  $G$-bundle   $E\to M$ with $G$ a compact semi-simple Lie
group, denote
by $\cG:=\Gamma(\Aut E)$ its gauge group with $\fg:=\Lie(\cG)$ the associated adjoint bundle
and by $\cA(E)$ its space of connections.
We address the classical problem of describing the absolute minima of the \emph{Yang-Mills
 functional}
\begin{eqnarray*}
\cS_{\rm YM}&\colon& \cA(E) \to \R^+\\
\cS_{\rm YM}(A)
&:=&\left\Vert F_A \right\Vert^2_{L^2(M)}
=\int_M \left\langle F_A\wedge*F_A\right\rangle_\frg
\end{eqnarray*}
i.e., solutions of the \emph{Yang-Mills equation}:
\begin{equation}
\label{eq: YM}
\de^*_A F_A=0.
\end{equation}

\subsection{Yang-Mills connections, $\rG_2$--instantons and the Chern-Simons action}
\label{sec: YM and CS}
The paradigmatic PDE for gauge theory in the presence of a $\rG_2$--structure is the \emph{$\rG_2$--instanton equation} \cites{Donaldson1998,Tian2000},
which can be formulated equivalently in terms of $\varphi$ or $\psi:=*\varphi$:
\begin{equation}
\label{eq: g2-instanton}
        F_A\wedge\psi =0\\
        \quad\Leftrightarrow\quad
        *F_A= F_A\wedge\varphi.
\end{equation}
This is the natural Euler-Lagrange equation for the  \emph{Chern-Simons action},
defined relatively to  a fixed reference connection $A_0\in\cA(E)$ by 
\begin{eqnarray*}
        \cS_{\rm CS}&\colon& \cA( E)\simeq A_0+\Omega^1(\frg) 
        \to \R\\
        \cS_{\rm CS}(A_0+a)
        &:=&\frac{1}{2}\int_M\tr 
        \left(\de_{A_0}a\wedge a
        +\frac{2}{3}a\wedge a\wedge a\right) 
        \wedge\psi
\end{eqnarray*}
with  $\cS_{\rm CS}(A_0)=0$.  Supposing the cocalibrated condition  $\de\psi=0$, the action is well-defined and    $\rG_2$--instantons are manifestly critical points. Its gradient is the \emph{Chern-Simons $1$--form} $\rho=\de\cS_{\rm CS}$, defined on vector fields $b\colon \cA( E)\to\Omega^1(\fg)$ by
\begin{equation}
\label{eq: CS 1-form}
\rho(b)_A=\int_M \tr \left(F_A\wedge b_A\right)\wedge\psi.
\end{equation}
and indeed  the solutions of  (\ref{eq: g2-instanton}) are precisely its zeroes (for a more detailed exposition, see \cite{SaEarp2014}.)  

Now, if the $\rG_2$--structure was closed, then by the Bianchi identity every solution of (\ref{eq: g2-instanton}) would automatically solve (\ref{eq: YM}). In other words, $\rG_2$--instantons would  manifestly be critical points of the Yang-Mills functional, somewhat in analogy to (anti-)selfdual connections
in dimension $4$ \cite{Donaldson1990}. This indeed was the starting point
of our predecessors in proposing gauge theory on $\rG_2$--manifolds. Since
then, such Yang-Mills  $\rG_2$--instantons have been constructed on Joyce manifolds \cite{Walpuski2013}, Bryant-Salamon manifolds \cite{Clarke2014}, associative fibrations
\cite{SaEarp2014}, asymptotically cylindrical  $\rG_2$--manifolds
\cites{SaEarp2009,SaEarp2015a} and their twisted connected sums  \cite{SaEarp2015b,Walpuski2015}. However, the implication (\ref{eq: g2-instanton}) $\Rightarrow$ (\ref{eq: YM}), that   $\rG_2$--instantons are Yang-Mills minima,  does not hold for generic  cocalibrated $\rG_2$--structures. In particular, a naive attempt to apply standard Chern-Weil arguments would depend
on a certain characteristic class in de Rham cohomology
which generically does not exist\footnote{The implication may still hold, in special circumstances, despite the absence of the de Rham class, see for instance the scalar torsion case in \cite{Harland2012}.}. This is unfortunate,
since every oriented spin $7$-manifold admits such a structure \cite{Crowley2014}, of which indeed many explicit examples are now known (see e.g. \cites{Conti2006,Agricola2010,Lotay2012,Freibert2013}
and references therein).  As we will show in 
Section  \ref{sec: top energy bounds}, a suitable version of the argument \emph{does}
hold for the natural cocalibrated  $\rG_2$--structure on contact Calabi-Yau
manifolds and in particular for Calabi-Yau links. Explicitly,  integrable $\rG_2$--instantons on holomorphic Sasakian vector bundles are solutions of the Yang-Mills equation  (\ref{eq: YM}), even though the   $\rG_2$--structure is not closed. 

\subsection{Sasakian vector bundles}
\label{sec: Sasakian vb}

Interpreting Sasakian structures as a setup for `transversely K\"ahler' geometry,
any compatible formulation of gauge theory requires a good notion of transversely
compatible structures on vector bundles. We adopt the lexicon proposed by
Biswas and Schumacher \cite[\S3.3]{Biswas2010}. 

Let $E\to M$ be a   $C^\infty$ complex vector bundle
over a smooth manifold and let
\begin{equation}
\label{eq: subbundle inclusion}
        L\subset TM_\C:=TM\otimes_{\R}\C
\end{equation}
be an integrable subbundle, i.e., closed under the Lie bracket. A \emph{partial
connection along $L$ }is a  $C^\infty$ differential operator
\begin{eqnarray*}
        \De_L:E &\longrightarrow& L^*\otimes E
\end{eqnarray*}
satisfying the Leibniz condition for $ f\in C^\infty(M)$ and $s\in\Gamma(E)$:
$$
\De_L(fs)=f\De_L(s)+q_L(\de f)
$$  
relative to the dual $q_L:T^*M_\C\to L^*$ of the inclusion (\ref{eq: subbundle inclusion}).
Since $L$ is integrable, $q_L$ induces a natural exterior derivative
$\hat \de : L^* \to\bigwedge^2 L^*$, hence an extension of $\De_L$ to $E$--valued sections of $L^*$
\begin{eqnarray*}
        \De_L:L^*\otimes E & \longrightarrow& ({\bigwedge}^2 L^*)\otimes
        E\\
        \De_L(t\otimes s) &:=& \hat\de(t)\otimes s-t\otimes \De_L(s).
\end{eqnarray*}
Denoting by $\fg$ the adjoint bundle of $E$, the \emph{curvature} of $\De_L$ is the $C^\infty(M)$--linear section
$$
F_{\De_L}:=\De_L^2 \in \Gamma \left( ({\bigwedge}^2 L^*)\otimes \fg \right)
$$
and $D_L$ is said to be \emph{flat} if $F_{\De_L}=0$. 
 
\begin{definition} 
\label{def: Sasakian vb}
 A \emph{Sasakian (vector) bundle} $\bE\to M$ over a Sasakian manifold  $\left(M, \theta, g, \Phi \right)$ with Reeb field $\xi$ is a pair $(E,\Dex)$ given
by a $C^\infty$ complex vector bundle $E$ over $M$ and a partial connection $\Dex:E\to\theta\otimes E$ along $\xi$.
\end{definition}
To be completely precise, we are applying the previous discussion to the
line subbundle  $L=\xi_\C:=N_\xi\otimes_\R\C\subset TM_\C$   spanned by $\xi$ over $\C$ (cf. Section \ref{sec: contact Calabi Yau}); it is
clear that any such $\Dex$ is flat. Moreover, the natural partial connection induced by $\Dex$ on $E^*$ gives natural definitions of $\bE^*$ and $\End({\bE})$
as Sasakian bundles, and we define a \emph{Hermitian structure} on $\bE$
as a smooth Hermitian structure $h$ on $E$ preserved by  $\Dex$, in the sense
that
$$
\de( h(s_1,s_2))|_{\xi_\C}=h(\Dex (s_1),s_2)+h(s_1,\Dex(s_2)).
$$ 
Clearly a Hermitian structure on  $\bE$ induces a Hermitian structure on  $\bE^*$ and on $\End({\bE})$.
A \emph{unitary connection}  on $(\bE,h)$ is a connection $A$ on $E$ 
such that $\de_A$ preserves $h$ in the usual sense.

Finally, in the notation of (\ref{eq: B^(p,q)}), we obtain a natural notion of transversely holomorphic  structures
over a  Sasakian manifold $M^{2n+1}$ relative to the  integrable `extended anti-holomorphic' $(n+1)$--dimensional foliation
\begin{equation}
\label{eq: hol distribution}
        \tilde B^{0,1}:=  B^{0,1} \oplus \xi_\C \subset TM_\C.
\end{equation}
\begin{definition}
\label{def: hol Sasakian vb}
A \emph{holomorphic (Sasakian)  bundle} $\cE\to M$ over a Sasakian manifold  $M $ with Reeb field $\xi$ is a pair $(\bE,\bar\partial)$ given by a Sasakian bundle $\bE=(E,\Dex)$ (cf. Definition \ref{def: Sasakian vb}) and
a flat partial connection $\bar\partial=\De_{\tilde B^{0,1}}$ such that $\bar\partial|_{\xi_\C}=\Dex$. 
\end{definition}  

An \emph{integrable connection} on $\cE=(\bE,\bar\partial)$ is connection $A$ on $E$
such that its induced partial connection along $\tilde B^{0,1}$ given by
$D_{\tilde B^{0,1}}:=\de_A|_{\tilde B^{0,1}}$ coincides
with $\bar\partial$.
We denote by $\cA(\cE)$ the subset of integrable connections inside $\cA(E)$.
We are now in position to extend the well-known concept of a  \emph{Chern connection}, mutually compatible with the holomorphic structure and the Hermitian
metric \cite[Proposition 2.1.56]{Donaldson1990}: 

\begin{proposition}[\cite{Biswas2010}, p.552]
\label{prop: Chern connection}
Let $(\cE,h)$ be a holomorphic Sasakian  bundle with Hermitian structure;
then there exists a unique unitary and integrable \emph{Chern connection} $A_h$ on $\cE$,
and 
$$
F_{A_h}\in  \Omega^{1,1}(\fg). 
$$
Moreover, the expression
\begin{equation}
        \det\left( \id_E+\frac{\bi}{2\pi} F_{A_h}\right)
        =:
        \sum_{j=0}^{n} c_j(\cE,h)
\end{equation}
defines closed \emph{Chern forms} $c_j(\cE,h) \in \Omega^{j,j}(M)$.  
\end{proposition}

 In a local holomorphic
trivialization $\tau$, i.e., such that $\cE$ is locally spanned by sections
in $\ker\delbar$, the Chern connection of $h$ is represented by the matrix of
$(1,0)$--forms $A^\tau_h=h^{-1}\partial h$ and its curvature has the form $F^\tau_{A_h}=\delbar(h^{-1}\partial
h)$. Moreover, it is clear from
(\ref{eq: hol distribution}) and Definition  \ref{def: hol Sasakian vb} that
 any other Hermitian structure $h'$ induces a Chern connection on $\cE$ satisfying $A_{h'}-A_h\in \Omega^{1,0}(\fg)$. 

\subsection{$\rG_2$--instantons and the Hermitian Yang-Mills condition}
\label{sec: G2 = tHYM}
A connection $A$ on a complex vector bundle over a K\"ahler manifold is \emph{Hermitian
Yang-Mills (HYM)} if 
$$
\hat F_A:=(F_A,\omega)= 0 
\qandq
F_A^{0,2}=0.
$$
This notion extends literally to Sasakian bundles $E \to M$, taking $\omega=d\theta\in\Omega^{1,1}(M)$ as the transverse K\"ahler form. Fixing a holomorphic structure on $E$, it is easy to check that compatible HYM connections
are exactly $\rG_2$--instantons:
\begin{lemma}
\label{lem: SD lifts to G2-instanton}
Let $\cE$ be a holomorphic Sasakian bundle over a $7$--dimensional closed contact Calabi-Yau manifold
$M$  endowed with its natural  $\rm G_2$--structure  (\ref{eq: cCY G2-structure}) . Then  a Chern connection $A$ on $\cE$ is HYM if, and only if, it is a $\rG_2-$instanton.
\begin{proof}
A Chern connection $A$ satisfies $F_{A}\in \Omega ^{1,1}\left( M \right)$
(Proposition \ref{prop: Chern connection}),
so taking account of the bidegree of the transverse holomorphic volume form
(cf. Definition \ref{def: contact CY mfd})
we have 
$%
F_{A}\wedge \epsilon =$ $F_{A}\wedge \bar{\epsilon }=0$ .  Therefore
\begin{displaymath}
        F_{A}\wedge \Im\epsilon 
        =\frac{1}{2\bi} F_{A}\wedge \left( \epsilon -\bar{\epsilon }\right) 
        =0.
\end{displaymath}
Now, taking the product with the $4$--form we have
$$
F_{A}\wedge\psi = \tfrac{1}{2} F_{A}\wedge \omega \wedge \omega=(cst.)\hat
F_A(*\theta),
$$
hence $A$ is a solution of (\ref{eq: g2-instanton}) if, and only if,
$\hat F_A=0$.
\end{proof}
\end{lemma} 

This result generalises the  well-known fact that HYM connections compatible
with a fixed holomorphic structure over a smooth Calabi-Yau $3$--fold pull back bijectively to $S^1$--invariant $\rG_2$--instantons over the product $CY^3\times S^1$ \cite[Proposition 8]{SaEarp2015a}.
Indeed, it is easy to deduce the corresponding claim for arbitrary circle
fibrations: 
\begin{corollary}
\label{cor: SD lifts to G2-instanton}
Let $X$ be a Calabi-Yau threefold, let $\pi:Y\to X$ be a Sasakian circle  fibration endowed with the natural $\rG_2$--structure (\ref{eq-GrayG2structure}), and let   $\cE:=\pi^*\cE_0\to Y$ be
the pullback from a holomorphic vector bundle  $\cE_0\to X$. Then $\cE$ is
a holomorphic Sasakian bundle, and  a Chern connection $A$ on $\cE_0$ is HYM if, and only if, $\pi^*A$
is a $\rG_2-$instanton on $\cE$.
\begin{proof}
The contact Calabi-Yau structure is trivially given by the global angular form $\theta\in\Omega^1(Y)$
and the pullbacks of the Calabi-Yau data from $X$, under $\pi$, with natural
Reeb field determined by $\theta(\xi)=1$, tangent to the $S^1$--action.  Then the
underlying complex vector bundle $\pi^*E_0$ is trivial along $\xi$ and we can
adopt $D_\xi=\de_\xi$  the trivial vertical connection, which is manifestly flat. This defines a Sasakian bundle structure (cf. Definition \ref{def: Sasakian vb}). Moreover, the $6$--dimensional distribution $B:=\ker\theta\subset TY$
maps under $\pi_*$ isomorphically to $TX$, which induces a natural bi-degree
decomposition $B=\bigoplus B^{i,j}$. It is immediate to check that the holomorphic structure  $\delbar_0$ on $\cE_0$ pulls back to a holomorphic structure $\delbar:=\pi^*\delbar_0$
on $\cE$.         
\end{proof}
\end{corollary} 
This gives a  correspondence
$$
\left\{\begin{tabular}{c}
$S^1$--invariant unitary  \\
connections on $\cE=\pi^*\cE_0$ \\
\end{tabular}\right\}
\overset{1-1}{\longleftrightarrow}
\left\{\begin{tabular}{c}
HYM Chern \\
connections  on $\cE_0$ \\
\end{tabular}\right\}
$$
which proves part \emph{(ii)} of Theorem \ref{thm: G2-intantons on K}. Notice
that the right-hand side is bijectively parametrised by stable holomorphic
structures on the underlying complex vector bundle  $E_0\to X$, by the Hitchin-Kobayashi
correspondence. 
At this point it is   sharply relevant to ask whether  a Sasakian
version of that correspondence may be obtained, via  a suitable notion of transverse stability. While Biswas and Schumacher do outline some progress in that direction \cite[\S3.4]{Biswas2010}, in the course of submission of the present article some fundamental progress has been achieved by Baraglia and Hekmati, who established a Hitchin-Kobayashi correspondence for transversely
K\"ahler foliated geometries \cite{Baraglia2018}.    

\begin{remark} 
\label{rem: cCY + G2-inst => YM}
The coincidence of $\rG_2$-instantons and HYM connections in the holomorphic Sasakian context of Lemma \ref{lem: SD lifts to G2-instanton} has as striking consequence on contact Calabi-Yau manifolds [cf. Proposition \ref{prop: G2-structure on cCY}]. Using the Bianchi identity, the Yang-Mills equation (\ref{eq: YM}) for any cocalibrated $\rG_2$-instanton $A$ [cf. (\ref{eq: g2-instanton})] is equivalently  rephrased as follows: \begin{eqnarray*}
        0=\de_A^*F_A=\de_A (F_A\wedge \varphi)=F_A\wedge \de\varphi 
        &\Leftrightarrow&
        F_A\in \ker \de\varphi\subset (\Omega^\bullet(M),\wedge).
\end{eqnarray*}
But in the cCY case, the natural  $\rG_2$-structure (\ref{eq: cCY G2-structure}) satisfies $\de\varphi=\omega^2$, and $F_A\in \ker \omega^2$ is precisely the HYM condition. Therefore integrable $\rG_2$-instantons are actually Yang-Mills critical points, even though the  $\rG_2$-structure is not closed. 
\end{remark} 

\subsection{Characteristic classes and topological energy bounds}
\label{sec: top energy bounds}
We will show that $7$-dimensional contact Calabi-Yau manifolds admit a naturally defined secondary characteristic class representing topological charge, which is another peculiar feature among $\rG_2$--structures
with torsion. From the perspective of gauge theory, this means that critical
points of the Chern-Simons functional indeed saturate the Yang-Mills energy,
just like in classical $4$-dimensional theory or more familiar torsion-free higher dimensional
models.

\begin{definition}
\label{def: charge k(A)}
Let $E\to M$ be a Sasakian  bundle (cf. Definition \ref{def: hol Sasakian vb}) over a $7$--dimensional closed contact Calabi-Yau manifold  (cf. Definition \ref{def: contact CY mfd}) with  $\rm G_2$--structure  (\ref{eq: cCY G2-structure}) given by Proposition \ref{prop: G2-structure on cCY}.
We define the \emph{charge} of a  connection $A\in\cA(E)$ by
\begin{equation}
\label{eq: k(E)}
        \kappa(A):=\int_M \tr F_A^2 \wedge\varphi
\end{equation}
\end{definition}

\begin{lemma}
\label{lemma: topological charge}
In the context of Proposition \ref{prop: Chern connection} and Definition \ref{def: charge k(A)}, the quantity (\ref{eq: k(E)}) assessed among Chern connections in $\cA(\cE)$ is independent 
of the Hermitian structure and it defines a \emph{topological charge} $\kappa(\cE)$.
\begin{proof}
Fix a Hermitian metric
$h$ on a holomorphic Sasakian bundle $\cE$ with Chern connection $A=A_h\in\cA(\cE)$; then any other Chern connection on $\cE$ has the form  $A'=A +b $, for some $b\in\Omega^{1,0}(\fg)$. We know from standard Chern-Weil
theory that$$
\tr F^2_{A+b}-\tr F^2_A = \de\left(\tr\eta\right)
$$
for some
$$
\eta=\eta(A,b)
:=F_A\wedge b
+\frac{1}{2}\de_{A}b\wedge b
+\frac{1}{3}b\wedge b\wedge b
\in\Omega^3(\fg).
$$
Since  by assumption $M$ is a closed manifold, the quantity $\kappa$ in (\ref{eq: k(E)}) is defined up to a term
given by Stokes' theorem after integration by parts:
\begin{equation}
\label{eq: Chern-Weil term}
        \int_M \tr\eta\wedge \de\varphi=\int_M \tr \left(F_{A}\wedge b
        +\frac{1}{2}\de_{A}b\wedge b+\frac{1}{3}b\wedge b\wedge b\right)\wedge \de\varphi.         
\end{equation}

Since $A$ is a Chern connection, Proposition   \ref{prop: Chern connection} specifies the bi-degree of 
$$
F_A\wedge b\in\Omega^{2,1}(\fg).
$$
Moreover, since $b$ is transversal to the characteristic foliation and $A$ is locally of type $(1,0)$, the term $d_Ab$ has basic type $(2,0)+(1,1)$, so
$$
\de_Ab\wedge b\in (\Omega^{3,0}\oplus\Omega ^{2,1})(\fg).
$$
Finally, one obviously has $b\wedge b\wedge b\in \Omega^{3,0} $.

Recall from Lemma \ref{lemma: dtheta is (1,1)}   and Proposition \ref{prop: G2-structure on cCY} that      $\de\varphi=(\de\theta)^2\in\Omega^{2,2}(M)$,
hence all three terms on the right-hand side  of (\ref{eq: Chern-Weil term}) vanish by excess in bi-degree.  

\end{proof}
\end{lemma}
Now, following a classical argument, on one hand we have the orthogonal decomposition of the Yang-Mills functional:%
\begin{equation}        \label{YM(A)}
        \cS_{\rm YM}\left( A\right) 
        = \Vert F_{A}\Vert ^{2}
        =\Vert F_7\Vert ^{2}
        +\Vert F_{14}\Vert ^{2}.
\end{equation}
On the other hand, applying the $\rG_2$--equivariant eigenspace decomposition  from Remark \ref{rem: decomposition} to integrable connections as in $(ii)$
of Lemma
\ref{lemma: topological charge}, a straightforward calculation relates the
topological charge to these components: %
\begin{displaymath}
        \kappa \left( \cE\right)
        =
        -2\left\Vert F_{7}\right\Vert^{2}+\left\Vert F_{14}\right\Vert ^{2}.
\end{displaymath}
Combining with (\ref{YM(A)}),
we can isolate the topological charge as a lower
bound of
the Yang-Mills energy among integrable connections:
\begin{equation}
\label{eq: topological energy}
        \cS_{\rm YM}|_{\cA(\cE)}(A)=-\frac{1}{2}\kappa \left( \cE\right)+\frac{3}{2}\Vert F_{14}\Vert ^{2} 
        =\kappa (\cE)+3\left\Vert F_{7}\right\Vert ^{2}.
\end{equation}%
Hence, if $\cS_{\rm YM} $ attains on $\cA(\cE)$ its absolute topological minimum, this occurs at a
connection whose curvature lies either in 
$\Omega_{7}^{2}$ or in 
$\Omega_{14}^{2}$. Moreover, since $\cS_{\rm YM}\geq0$, the sign of  $\kappa(\cE)$ obstructs the existence of one type or the other, so we fix $\kappa(\cE)\geq0$,
compatibly with the existence of our $\rG_2-$instantons  (\ref{eq: g2-instanton})
with $ F_7=0$, i.e., such that $\cS_{\rm YM}(A)=\kappa(\cE)$. Together with Remark \ref{rem: cCY + G2-inst => YM}, this proves Theorem \ref{thm: G2-inst are YM minima}.

In summary, under the natural cocalibrated  $\rG_2$--structure of a contact Calabi-Yau $7$--manifold,   $\rG_2$--instantons (\ref{eq: g2-instanton}) over  
are Yang-Mills critical points and indeed \emph{absolute minima} among compatible connections.

\subsection{Example:   $\rG_2-$instantons on pullback bundles}
\label{sec: example pullback}

Motivated by Corollary \ref{cor: SD lifts to G2-instanton} in the previous section, let us explore the simplest model case for gauge theory on a $7$--dimensional
CY
link.
Let $\pi:Y=K_f \to V$ be a CY link, which fibres nontrivially by circles over the smooth  $3$--fold $V=(f)\subset\P^4$, and consider the holomorphic
Sasakian bundle given by pullback $\cE:=\pi^*\cE_0\to K_f$
 of a holomorphic bundle over $V$.
We would like to describe the explicit  local form of the constraint imposed on a Chern connection
 $\bA\in \cA(\cE)$ by the   $\rG_2-$instanton equation (\ref{eq: g2-instanton}). 

Over a trivialising
neighbourhood of $K_f $ as a circle fibration, i.e. an open set $ U\subset V$ such that $K_f \supset\ \pi^{-1}(U)\simeq
S^1 \times U$, given points $y\in  \pi^{-1}(U)$ and $x=\pi(y)\in U$, an
arbitrary integrable connection
$\mathbf{A}$ on  $\cE$ can be written as 
\begin{displaymath}
        \mathbf{A}(y)\overset{\loc}{=}
        \pi^*A_t(x)+ \sigma(x,t)\theta
\end{displaymath}
where $\left\{A_t\right\}_{t\in S^{1}}$ is a family of integrable connections on $\cE_0$
and   $\sigma\in\Omega^0(K,\fg)$, where $\fg:=\pi^*\mathfrak{g}_{\cE_0}$
is the corresponding adjoint bundle of $\cE$. Let us denote this fact
informally by
$$
\mathbf{A}=A_t+\sigma\theta.
$$ The curvature of $\mathbf{A}$ is the gauge-covariant global
$2$--form
\begin{displaymath}
        F_\mathbf{A} = F_{A_t} 
        \oplus\left(
        \de_{A_t}\sigma-\frac{\partial A_t}{\partial t} 
        \right)
        \wedge \theta\in\Omega^{2}(K,\fg).
\end{displaymath} 
and, replacing that expression in the $\rG_2$--instanton equation  (\ref{eq: g2-instanton}), one obtains
in particular$$
\hat F_{A_t}(*\theta)=F_{A_t}\wedge\omega^2=0.
$$
This is exactly the HYM condition on each $A_t$.
On the other hand, by Theorem \ref{thm: G2-inst are YM minima}, if  $\mathbf{A}$ is an integrable 
 $\rG_2-$instanton, then it minimises
the Yang-Mills functional (\ref{YM(A)}).
This implies locally
\begin{displaymath}
        \left(
        \de_{A_t}\sigma-\frac{\partial A_t}{\partial t} 
        \right)
        \wedge \theta
        =0,
\end{displaymath}
since otherwise  the pullback component $A_t$ alone would violate the minimum
topological energy (\ref{eq: topological energy}):
\begin{displaymath}
\cS_{\rm YM}(A_t) =\left\Vert F_{A_t}\right\Vert^2\geq\left\Vert F_{\mathbf{A}}\right\Vert^2=\cS_{\rm YM}(\mathbf{A})=\kappa(\cE).
\end{displaymath}
Moreover, if the family $A_t\equiv A_{t_0}$ is constant,  i.e., $S^1$--invariant,
then $\de_{A_{t_0}}\sigma=0$
 implies
   $\sigma\equiv0$, since by assumption $\cE$ is indecomposable and
therefore does not admit nonzero parallel sections, and so  $\mathbf{A}$ is indeed a pullback.
If the moduli space 
$\hat\cM$ of HYM connections on the base $V$ is discrete, then by continuity the family $\left\{A_t\right\}$
is contained in a gauge orbit. This concludes the proof of Theorem \ref{thm: G2-intantons on K}.
\begin{remark}  Let $\pi:Y=K_f \to V$ be the Fermat quintic link (cf. Example \ref{exa: Fermat quintic}). Since the moduli space of stable holomorphic bundles on a Fermat quintic Calabi-Yau $3$--fold $V$ is known to be discrete, we infer that  $S^1$--invariant  $\rG_2-$instantons should be counted in some sense by the Donaldson-Thomas invariant
of $V$, which is deformation-invariant because $h^{0,2}(V)=0$ 
\cite[Definition 3.34]{Thomas2000}. Thus we envisage a
`conservation of number' property for  $S^1$--invariant  $\rG_2-$instantons
over such Fermat quintic links, to be made precise in upcoming work. \end{remark}

\section*{Afterword: Atiyah's conjecture and singular   $\rm G_2$--metrics  }

Atiyah predicted that the Casson invariant $\lambda(\Sigma)$ of a homology sphere which is the link of a normal complete intersection singularity
equals $\tfrac{1}{8}\sigma(F)$, where $F$ is the Milnor fibre. This
was verified for Brieskorn spheres by Fintushel and Stern \cite{Fintushel1990},
and Neumann and Wahl \cite{Neumann1990} inductively use that fact to confirm the conjecture for weighted homogeneous surface singularities and for links of hypersurfaces of the form $f(x,y)+z^n=0$,
among others. Their theorem suggests a general relation between the Floer homology (or at least the Casson invariant) of a 
link in $\C^3$ and the signature of $F$. Arnold and  Floer \cite{Arnold1995} suggested higher-dimensional analogues, which  would require extra structure on the links (e.g. CR or contact structure) and Milnor fibre (e.g. symplectic structure). 

In our context, Chern-Simons theory (\ref{eq: CS 1-form}) suggests thinking of $\rG_2$--instantons as $7$--dimensional analogues of flat connections.
Applying the above  intuition to
the holomorphic Casson invariant of  R. Thomas over a CY $3$--fold
base \cite{Thomas2000}, we  wonder whether  a version of Atiyah's conjecture may hold for CY links. 

Finally, from the perspective of M-theory,  examples of compact  $\rm G_2$--metrics with prescribed singularities might be within reach,  starting from some suitably singular CY link and taking adiabatic limits
on the circle fibres near an orbifold singularity. As we have shown, meaningful Yang-Mills theory results may be established on such spaces, even though
the $\rG_2$--structure has some torsion. 

\appendix
\section{Algorithm for Steenbrink's signature theorem}
\label{app: algorithm}

As discussed in Section \ref{sec: Formula for nu}, Steenbrink's method for the signature of a compactified affine variety depends solely on computing the nonzero signature  $(\mu_{+},\mu_{-})$. This requires  an explicit basis of the Milnor algebra, which several computational tools provide. We use the following code in \textsc{Singular} \cite{Greuel2001}.

We first compute the numbers (\ref{eq-array Steenbrink}), for a given polynomial $f$ of degree $d$ and a list of weights $W$, and arrange them
into list $L$: \\

\noindent
\tt{proc Signature(poly f, list W, int d)}\\
\tt{ \{ }\\
\tt{ring A = 0, (a,b,c,d,e), lp;}\\
\tt{list L;}\\
\tt{int s;}\\
\tt{ideal J = jacob(f);}\\
\tt{J = groebner(J);}\\
\tt{ideal K = kbase(J);}\\
\tt{s=size(K);}\\
\tt{for ( int j=1; j <= s; j++ )}\\
\tt{ \{ L[j]=(1+leadexp(K[j])[1])*(W[1]/d) }\\
\tt{+ (1+leadexp(K[j])[2])*(W[2]/d)+ (1+leadexp(K[j])[3])*(W[3]/d)}\\
\tt{+ (1+leadexp(K[j])[4])*(W[4]/d)+(1+leadexp(K[j])[5])*(W[5]/d); \} }\\
\tt{return(L);}\\
\tt{write("list.txt", L);}\\ 
\tt{ \} }\\
\normalfont
 Then we use  \textsc{Mathematica}  code to compute the $\nu$ invariant from the list $L$:\\

$\nu$=\texttt{ Mod[Length[L]+1-3*(Length[Select[Select[L, \# \textbackslash[NotElement] Integers \&], 
\\Mod[IntegerPart[\#], 2] == 0 \&]]}
\tt{-Length[Select[Select[L, 
\\ \# \textbackslash[NotElement] Integers \&], Mod[IntegerPart[\#], 2] == 1 \&]]),48]}\\


\bibliography{Bibliografia-2018-06}

\begin{bibdiv}
\begin{biblist}

\bib{Agricola2010}{article}{
      author={Agricola, Ilka},
      author={Friedrich, Thomas},
       title={3-{S}asakian manifolds in dimension seven, their spinors and
  {$G_2$}-structures},
        date={2010},
        ISSN={0393-0440},
     journal={J. Geom. Phys.},
      volume={60},
      number={2},
       pages={326\ndash 332},
      review={\MR{2587396}},
}

\bib{Acharya1997}{article}{
      author={Acharya, B.~S.},
      author={O'Loughlin, M.},
      author={Spence, B.},
       title={Higher-dimensional analogues of {D}onaldson-{W}itten theory},
        date={1997},
        ISSN={0550-3213},
     journal={Nuclear Phys. B},
      volume={503},
      number={3},
       pages={657\ndash 674},
      review={\MR{1482739}},
}

\bib{Arnold1995}{incollection}{
      author={Arnold, V.~I.},
       title={Some remarks on symplectic monodromy of {M}ilnor fibrations},
        date={1995},
   booktitle={The {F}loer memorial volume},
      series={Progr. Math.},
      volume={133},
   publisher={Birkh\"auser, Basel},
       pages={99\ndash 103},
      review={\MR{1362824 (96m:32043)}},
}

\bib{Boyer2008}{book}{
      author={Boyer, Charles~P.},
      author={Galicki, Krzysztof},
       title={Sasakian geometry},
      series={Oxford Mathematical Monographs},
   publisher={Oxford University Press, Oxford},
        date={2008},
        ISBN={978-0-19-856495-9},
      review={\MR{2382957 (2009c:53058)}},
}

\bib{Baraglia2018}{article}{
      author={Baraglia, D.},
      author={Hekmati, P.},
       title={A foliated hitchin-kobayashi correspondence},
        date={2018},
     journal={arXiv:1802.09699 [math.DG]},
}

\bib{Bryant1987}{article}{
      author={Bryant, R.~L.},
       title={Metrics with exceptional holonomy},
        date={1987},
        ISSN={0003-486X},
     journal={Ann. of Math. (2)},
      volume={126},
      number={3},
       pages={525\ndash 576},
      review={\MR{916718 (89b:53084)}},
}

\bib{Biswas2010}{article}{
      author={Biswas, I.},
      author={Schumacher, G.},
       title={Vector bundles on {S}asakian manifolds},
        date={2010},
     journal={Adv. Theor. Math. Phys.},
      volume={14},
      number={2},
       pages={541\ndash 562},
}

\bib{Corrigan1983}{article}{
      author={Corrigan, E.},
      author={Devchand, C.},
      author={Fairlie, D.~B.},
      author={Nuyts, J.},
       title={First-order equations for gauge fields in spaces of dimension
  greater than four},
        date={1983},
        ISSN={0550-3213},
     journal={Nuclear Phys. B},
      volume={214},
      number={3},
       pages={452\ndash 464},
      review={\MR{698892 (84i:81058)}},
}

\bib{Crowley2015}{article}{
      author={Crowley, Diarmuid},
      author={Goette, Sebastian},
      author={Nordstr{\"o}m, Johannes},
       title={{An analytic invariant of $\rG_2$--manifolds}},
        date={2015},
       pages={1\ndash 26},
      eprint={arXiv:1505.02734v1},
}

\bib{Corti2013}{article}{
      author={Corti, A.},
      author={Haskins, M.},
      author={Nordstr{\"o}m, J.},
      author={Pacini, T.},
       title={{Asymptotically cylindrical Calabi--Yau $3$--folds from weak Fano
  $3$--folds}},
        date={2013},
     journal={Geom. Topol.},
      volume={17},
      number={4},
       pages={1955\ndash 2059},
}

\bib{Corti2015}{article}{
      author={Corti, A.},
      author={Haskins, M.},
      author={Nordstr{\"o}m, J.},
      author={Pacini, T.},
       title={{$\rG_2$--manifolds and associative submanifolds via semi-Fano
  $3$--folds}},
        date={2015},
     journal={Duke Math. J.},
      volume={164},
      number={10},
       pages={1971\ndash 2092},
      eprint={arxiv:1207.4470v2},
}

\bib{Clarke2014}{article}{
      author={Clarke, Andrew},
       title={Instantons on the exceptional holonomy manifolds of {B}ryant and
  {S}alamon},
        date={2014},
     journal={Journal of Geometry and Physics},
      volume={82},
       pages={84\ndash 97},
}

\bib{Candelas1990}{article}{
      author={Candelas, P.},
      author={Lynker, M.},
      author={Schimmrigk, R.},
       title={Calabi-{Y}au manifolds in weighted {${\bf P}_4$}},
        date={1990},
        ISSN={0550-3213},
     journal={Nuclear Phys. B},
      volume={341},
      number={2},
       pages={383\ndash 402},
      review={\MR{1067295}},
}

\bib{Crowley2014}{article}{
      author={Crowley, Diarmuid},
      author={Nordstr{\"o}m, Johannes},
       title={Exotic {$G_2$}-manifolds},
        date={2014},
     journal={arXiv:1411.0656 [math.AG]},
}

\bib{Crowley2015b}{article}{
      author={Crowley, Diarmuid},
      author={Nordstr{\"o}m, Johannes},
       title={New invariants of {$\rG_2$}--structures},
        date={2015},
     journal={Geometry \& Topology},
      volume={19},
      number={5},
       pages={2949\ndash 2992},
      review={\MR{3416118}},
}

\bib{Conti2006}{article}{
      author={Conti, D.},
      author={Salamon, S.},
       title={Reduced holonomy, hypersurfaces and extensions},
        date={2006},
     journal={Int. J. Geom. Methods Mod. Phys.},
      volume={3},
}

\bib{Dimca1992}{book}{
      author={Dimca, Alexandru},
       title={Singularities and topology of hypersurfaces},
      series={Universitext},
   publisher={Springer-Verlag, New York},
        date={1992},
        ISBN={0-387-97709-0},
      review={\MR{1194180}},
}

\bib{Donaldson2002}{book}{
      author={Donaldson, S.~K.},
       title={Floer homology groups in {Y}ang--{M}ills theory},
      series={Cambridge Tracts in Mathematics},
   publisher={Cambridge University Press},
     address={Cambridge},
        date={2002},
      volume={147},
        ISBN={0-521-80803-0},
        note={With the assistance of M. Furuta and D. Kotschick},
      review={\MR{1883043 (2002k:57078)}},
}

\bib{Donaldson1990}{article}{
      author={Donaldson, S.~K.},
       title={Polynomial invariants for smooth four--manifolds},
        date={1990},
     journal={Topology},
      volume={29},
      number={3},
       pages={257\ndash 315},
}

\bib{Donaldson1998}{incollection}{
      author={Donaldson, S.~K.},
      author={Thomas, R.~P.},
       title={Gauge theory in higher dimensions},
        date={1998},
   booktitle={The geometric universe ({O}xford, 1996)},
   publisher={Oxford Univ. Press},
     address={Oxford},
       pages={31\ndash 47},
         url={http://www.ma.ic.ac.uk/~rpwt/skd.pdf},
      review={\MR{MR1634503 (2000a:57085)}},
}

\bib{Fernandez1982}{article}{
      author={Fern{\'a}ndez, M.},
      author={Gray, A.},
       title={Riemannian manifolds with structure group $\rm{G}_2$},
        date={1982},
        ISSN={0003-4622},
     journal={Annali di Matematica Pura ed Applicata},
      volume={132},
      number={1},
       pages={19\ndash 45 (1983)},
      review={\MR{696037 (84e:53056)}},
}

\bib{Freibert2013}{article}{
      author={Freibert, Marco},
       title={{Cocalibrated $\rG_2$--structures on products of four- and
  three-dimensional {L}ie groups}},
        date={2013},
     journal={Differential Geom. Appl.},
      volume={31},
      number={3},
       pages={349\ndash 373},
}

\bib{Fintushel1990}{article}{
      author={Fintushel, Ronald},
      author={Stern, Ronald~J.},
       title={Instanton homology of {S}eifert fibred homology three spheres},
        date={1990},
        ISSN={0024-6115},
     journal={Proc. London Math. Soc. (3)},
      volume={61},
      number={1},
       pages={109\ndash 137},
      review={\MR{1051101 (91k:57029)}},
}

\bib{Greuel2001}{techreport}{
      author={Greuel, G.-M.},
      author={Pfister, G.},
      author={Sch\"onemann, H.},
       title={{\sc Singular} 2.0},
        type={{A Computer Algebra System for Polynomial Computations}},
 institution={Centre for Computer Algebra},
     address={University of Kaiserslautern},
        date={2001},
        note={{\tt http://www.singular.uni-kl.de}},
}

\bib{Gray1969}{article}{
      author={Gray, A.},
       title={Vector cross products on manifolds},
        date={1969},
        ISSN={0002-9947},
     journal={Trans. Amer. Math. Soc.},
      volume={141},
       pages={465\ndash 504},
      review={\MR{0243469 (39 \#4790)}},
}

\bib{Harvey1982}{article}{
      author={Harvey, F.~R.},
      author={Lawson, H.~B., Jr.},
       title={Calibrated geometries},
        date={1982},
        ISSN={0001-5962},
     journal={Acta Math.},
      volume={148},
       pages={47\ndash 157},
      review={\MR{MR666108 (85i:53058)}},
}

\bib{Harvey1999}{article}{
      author={Harvey, Jeffrey~A.},
      author={Moore, Gregory~W.},
       title={{Superpotentials and membrane instantons}},
        date={1999},
      eprint={hep-th/9907026},
}

\bib{Harland2012}{article}{
      author={Harland, D.},
      author={N{\"o}lle, C.},
       title={Instantons and killing spinors},
        date={2012-03},
     journal={Instantons and Killing spinors},
      volume={82},
}

\bib{Habib2015}{article}{
      author={Habib, Georges},
      author={Vezzoni, Luigi},
       title={Some remarks on {C}alabi-{Y}au and hyper-{K}\"ahler foliations},
        date={2015},
        ISSN={0926-2245},
     journal={Differential Geom. Appl.},
      volume={41},
       pages={12\ndash 32},
      review={\MR{3353736}},
}

\bib{Joyce2000}{book}{
      author={Joyce, D.~D.},
       title={Compact manifolds with special holonomy},
      series={Oxford Mathematical Monographs},
   publisher={Oxford University Press},
     address={Oxford},
        date={2000},
        ISBN={0-19-850601-5},
      review={\MR{1787733 (2001k:53093)}},
}

\bib{Lotay2012}{article}{
      author={Lotay, J.~D.},
       title={Associative submanifolds of the $7$-sphere},
        date={2012},
     journal={Proceedings of the London Mathematical Society},
      volume={105},
      number={6},
       pages={1183\ndash 1214},
}

\bib{Milnor1969}{misc}{
      author={Milnor, John},
       title={{Singular points of complex hypersurfaces}},
   publisher={Princeton University Press},
        date={1969},
}

\bib{Milnor1970}{article}{
      author={Milnor, John},
      author={Orlik, Peter},
       title={Isolated singularities defined by weighted homogeneous
  polynomials},
        date={1970},
        ISSN={0040-9383},
     journal={Topology},
      volume={9},
       pages={385\ndash 393},
      review={\MR{0293680}},
}

\bib{Neumann1990}{article}{
      author={Neumann, Walter},
      author={Wahl, Jonathan},
       title={Casson invariant of links of singularities},
        date={1990},
        ISSN={0010-2571},
     journal={Comment. Math. Helv.},
      volume={65},
      number={1},
       pages={58\ndash 78},
      review={\MR{1036128 (91c:57022)}},
}

\bib{Portilla2020}{thesis}{
      author={Portilla, Luis},
       title={to appear},
        type={Ph.D. Thesis},
        date={2020},
}

\bib{Reinhart1959}{article}{
      author={Reinhart, B.},
       title={Foliated manifolds with bundle-like metrics},
        date={1959},
     journal={Annals of Mathematics},
      volume={69},
       pages={119\ndash 132},
}

\bib{Salamon1989}{book}{
      author={Salamon, S.},
       title={Riemannian geometry and holonomy groups},
      series={Pitman Research Notes in Mathematics Series},
   publisher={Longman Scientific \& Technical},
     address={Harlow},
        date={1989},
      volume={201},
        ISBN={0-582-01767-X},
      review={\MR{1004008 (90g:53058)}},
}

\bib{SaEarp2009}{thesis}{
      author={S{\'a}~Earp, H.~N.},
       title={Instantons on $\rm{G}_2-$manifolds},
        type={Ph.D. Thesis},
        date={2009},
}

\bib{SaEarp2014}{article}{
      author={S{\'a}~Earp, H.~N.},
       title={{Generalised Chern-Simons theory and $\rm G_2-$instantons over
  associative fibrations}},
        date={2014},
     journal={SIGMA - Symmetry, Integrability and Geometry: Methods and
  Applications},
      volume={10},
      number={083},
      eprint={arXiv:1401.5462v2 [math.DG]},
}

\bib{SaEarp2015a}{article}{
      author={S{\'a}~Earp, H.~N.},
       title={$\rm{G}_2-$instantons over asymptotically cylindrical manifolds},
        date={2015},
     journal={Geometry \& Topology},
      volume={19},
      number={1},
       pages={61\ndash 111},
}

\bib{SaEarp2015b}{article}{
      author={S{\'a}~Earp, H.~N.},
      author={Walpuski, T.},
       title={$\rm{G}_2-$instantons over twisted connected sums},
        date={2015},
     journal={Geometry \& Topology},
      volume={19},
      number={3},
       pages={1263\ndash 1285},
}

\bib{Sasaki1962}{article}{
      author={Sasaki, S.},
      author={Hatakeyama, Y.},
       title={On differentiable manifolds with contact metric structures},
        date={1962},
        ISSN={0025-5645},
     journal={J. Math. Soc. Japan},
      volume={14},
       pages={249\ndash 271},
      review={\MR{0141045 (25 \#4458)}},
}

\bib{Steenbrink1977}{article}{
      author={Steenbrink, Joseph},
       title={Intersection form for quasi-homogeneous singularities},
        date={1977},
        ISSN={0010-437X},
     journal={Compositio Math.},
      volume={34},
      number={2},
       pages={211\ndash 223},
      review={\MR{0453735}},
}

\bib{Thomas2000}{article}{
      author={Thomas, Richard},
       title={{A holomorphic {C}asson invariant for {C}alabi-{Y}au $3$--folds,
  and bundles on $K3$ fibrations}},
        date={2000},
     journal={Jour. Diff. Geom.},
      volume={54},
      number={2},
       pages={367\ndash 438},
}

\bib{Tian2000}{article}{
      author={Tian, G.},
       title={Gauge theory and calibrated geometry. {I}},
        date={2000},
        ISSN={0003-486X},
     journal={Ann. of Math. (2)},
      volume={151},
      number={1},
       pages={193\ndash 268},
      review={\MR{MR1745014 (2000m:53074)}},
}

\bib{Tomassini2008}{article}{
      author={Tomassini, Adriano},
      author={Vezzoni, Luigi},
       title={Contact {C}alabi-{Y}au manifolds and special {L}egendrian
  submanifolds},
        date={2008},
        ISSN={0030-6126},
     journal={Osaka J. Math.},
      volume={45},
      number={1},
       pages={127\ndash 147},
         url={http://projecteuclid.org/euclid.ojm/1205503561},
      review={\MR{2416653}},
}

\bib{Walpuski2013}{article}{
      author={Walpuski, T.},
       title={$\rm{G}_2-$instantons on generalised {K}ummer constructions},
        date={2013},
     journal={Geometry {\&} Topology},
      volume={17},
      number={4},
       pages={2345–\ndash 2388},
}

\bib{Walpuski2015}{article}{
      author={Walpuski, Thomas},
       title={{$\rG_2$--instantons over twisted connected sums : an example}},
        date={2015},
      eprint={arXiv:1505.01080v1},
}

\bib{Witten1988}{article}{
      author={Witten, Edward},
       title={Topological quantum field theory},
        date={1988},
        ISSN={0010-3616},
     journal={Comm. Math. Phys.},
      volume={117},
      number={3},
       pages={353\ndash 386},
         url={http://projecteuclid.org/euclid.cmp/1104161738},
      review={\MR{953828}},
}

\end{biblist}
\end{bibdiv}

\end{document}